\documentclass[a4paper,11pt]{article}
\usepackage{fullpage}
\usepackage[latin1]{inputenc}
\usepackage{graphicx}
\usepackage{indentfirst}
\usepackage[T1]{fontenc}
\usepackage{amsmath}
\usepackage{amssymb}
\usepackage{amsfonts}
\usepackage{latexsym}
\usepackage{graphicx}
\usepackage{fancybox}
\usepackage{epsfig}
\usepackage{psfig}
\usepackage{float}

\newcommand{\dsp}{\displaystyle}
\newcommand{\barr}{\begin{array}}
\newcommand{\earr}{\end{array}}
\newcommand{\Div}{\mbox{\rm div}}

\newtheorem{remark}{Remark}

\begin{document}
\titlepage
\title{A time domain method for modeling viscoacoustic wave propagation}
\author{ Jean-Philippe Groby\thanks{Laboratoire de Mécanique de d'Acoustique,
  Marseille, FRANCE, ({\tt groby@lma.cnrs-mrs.fr})}
 \and Chrysoula
Tsogka\thanks{Mathematics Department, Stanford University, USA, ({\tt
    tsogka@math.Stanford.EDU})}}
\date{\today}
\maketitle
\begin{abstract}
In many applications, and in particular in seismology, realistic
propagation media disperse and attenuate waves. This dissipative
behavior can be taken into account by using a viscoacoustic
propagation model, which incorporates a complex and 
frequency-dependent viscoacoustic modulus in the constitutive relation. The
main difficulty then lies in finding an efficient way to
discretize the constitutive equation as it becomes a convolution
integral in the time domain. To overcome this difficulty the usual
approach consists in approximating the viscoacoustic modulus by a
low-order rational function of frequency. We use here such an
approximation and show how it can be incorporated in the
velocity-pressure formulation for viscoacoustic waves. This
formulation is coupled with the fictitious domain method which
permit us to model efficiently diffraction by objects of
complicated geometry and with the Perfectly Matched Layer Model
which allows us to model wave propagation in unbounded domains.
The space discretization of the problem is based on a mixed finite
element method and for the discretization in time a 2nd order
centered finite difference scheme is employed. Several numerical
examples illustrate the efficiency of the method.
\end{abstract}

\section{Introduction}
Real media attenuate and disperse propagating waves \cite{Fut62}.
Our aim in this paper is to develop a numerical method to model
such dissipative phenomena (dispersion plus attenuation) in the
time domain. To do so we consider the linear viscoacoustic
equation which is a convolution in the time domain, the
viscoacoustic modulus being frequency dependent. Therefore,
incorporating any arbitrary dissipation law in time-domain methods
is in general computationally intense. The usual way to overcome
this difficulty is to approximate the viscoacoustic modulus by a
low-order rational function \cite{DM84, EK87, CKK88a, BRS95}. This
leads to replacing the convolution integral by a set of variables,
usually referred to as memory variables, which satisfy simple
differential equations that can be easily discretized in the time
domain.

Several methods have been proposed in the literature for
incorporating realistic attenuation laws (e.g.
frequency-independent or weakly frequency-dependent viscoacoustic
modulus) into time-domain methods \cite{DM84, EK87, BRS95, CKK88a,
CKK88b, CKK88c}. We focus our attention in this paper on the
methods proposed by  Day and Minster (1984), Emmerich and Korn
(1987), and Blanch, Robertson and Symes (1995). All three methods
use some approximation of the viscoacoustic modulus by a low-order
rational function. The first approach is based on the standard
Pad\'e approximation. The coefficients of the rational
approximation are thus in principle known analytically. Numerical
results obtained using this method show that the approximation is
poor and the method provides satisfactory results only for
relatively short (in terms of the wavelength) propagation paths.
The second approach is based on the rheological model of the
generalized Maxwell body, which gives a physical meaning to the
coefficients of the rational approximation. They are interpreted
as the relaxation frequencies and  weight factors of the classical
Maxwell bodies, which form the generalized Maxwell body. This
method provides good numerical results for long propagation paths,
but some parameters, namely the relaxation frequencies are
semi-empirically determined. Finally, the third method is based on
the observation that for the frequency-independent case and for
weakly attenuating materials the weight factors are only slowly
varying and can be approximated by a constant. This method
provides good numerical results, but also involves a
semi-empirical choice of a parameter.

Although, the previous methods give satisfactory results in the
case of weakly-attenuating materials they fail in media with large
attenuation. This case was considered in a recent paper
\cite{Asv2003}, where the authors propose an analytic method for
computing the best (optimal) rational approximation for the
frequency independent case. They also propose a generalization of
the algorithm presented in \cite{EK87} which leads to very good
results in the case of highly attenuating media and a frequency-
dependent viscoacoustic modulus.

After a brief overview of the basic theory describing wave propagation
in viscoacoustic media (section \ref{s1}), we describe in section \ref{s2}
the approximations proposed in \cite{DM84}, \cite{EK87} and
\cite{BRS95}.

Considering long propagation paths, we test the performance of the
different approximations and find that the best method, using the
smaller number of unknowns while providing satisfactory numerical
results and involving the least number of empirically determined
values, is the one proposed by Emmerich and Korn (1987). We thus
chose this method for approximating the viscoacoustic modulus.
Note that a slight variation of the method proposed in \cite{EK87}
is used here, based on a different way of distributing the
relaxation frequencies in the bandwidth of the incident pulse.

In section \ref{s3} we incorporate this approximation in the
velocity-pressure formulation for viscoacoustic waves.  Our choice
of using the first-order-in-time system of equations, instead of
the more classical second-order one, is motivated by the use of
the fictitious domain method and the perfectly matched absorbing
layer technique. In \cite{BEJ03} the authors proposed a similar
approach using the mixed velocity-stress formulation for modeling
wave propagation in viscoelastic media.

The fictitious domain method (also called the domain embedding
method) has been developed for solving problems involving complex
geometries \cite{B73, GPP94a, GPP94b, GG95, GK98}, and, in
particular, for wave propagation problems \cite{CJM97, Sylv, Leil,
BJT4}. In the framework of seismic wave propagation we apply this
method to model the boundary condition on the surface of the earth
(section \ref{s7}).  Its main feature is extending the solution to
a domain with simple shape, independent of the complex geometry,
and to impose the boundary conditions with the introduction of a
Lagrange multiplier. Thus, the solution is determined by two types
of unknowns, the extended unknowns, defined in the enlarged simple
shape domain and the auxiliary variable, supported on the boundary
of complex geometry. The main advantage is that the mesh for
computing the extended functions can  now be chosen independently
of the geometry of the boundary.

The  Perfectly Matched Layers (PML) technique was introduced by
B\'erenger \cite{ber94, ber96} for  Maxwell's equations and is now
the most widely-used method for the simulation of electromagnetic
waves in unbounded domains (cf. \cite{zhao,
  teix2, Petr}). It has also been extended to the case of anisotropic
acoustic waves \cite{BFP}, isotropic \cite{hastings} and anisotropic
elastic waves \cite{CT, BFP}. This technique consists in designing an
absorbing layer, called a perfectly matched layer (PML), that has
the property of generating {\it no reflection}
at the interface between the free medium and the artificial
absorbing medium. This property allows the use of a very high
damping parameter inside the layer, and consequently of a small
layer width, while achieving a near-perfect absorption of the
waves.  We apply here the PML model in the case of viscoacoustic
waves (section \ref{s6}).

Another advantage of the first-order formulation over the second
order one, is that it is easier to implement in heterogeneous
media, since it does not require an approximation of spatial
derivatives of the physical parameters. To discretize this
formulation in space we use a mixed finite-element method which is
a modification of the method proposed in \cite{BJT1-SIAM}. More
precisely, in \cite{BJT1-SIAM} the authors designed new mixed
finite elements, the so-called $Q_{k+1}^{div}-Q^{k}$ elements,
inspired  by Néd\'elec's second family \cite{Ned},  which are
compatible with mass lumping, and therefore allow to construct an
explicit scheme in time. A non-standard  convergence analysis of
the $Q_{k+1}^{div}-Q^{k}$ elements was carried out in
\cite{BJT1-SIAM}. However, numerical results obtained recently
(cf. \cite{BRT}) show that, when coupled with the fictitious
domain method, these elements do not provide satisfactory results.
This is why we use here instead the $Q_{k+1}^{div}-P^{k+1}$
elements for which convergence of the fictitious domain method was
obtained \cite{BRT}.

To show the efficiency and robustness of the method we present in
section \ref{s8} several numerical results. In particular,
numerical and analytical results are compared and  good agreement
is obtained between the two.
\section{Viscoacoustic wave propagation}{\label{s1}}
In an isotropic viscoacoustic medium occupying a domain $\Omega
\in \mathbb{R} ^{d}$, $d=1,2,3$, the relation between the pressure
$p(\omega)=p(\mathbf{x},\omega)$ and the displacement
$\mathbf{u}(\omega)=\mathbf{u}(\mathbf{x},\omega)$  in the
frequency domain is,
\begin{equation}
\dsp p(\omega)=\mu(\omega) \Div \mathbf{u}(\omega).
\label{relationpufreq}
\end{equation}
Here, $\mu(\omega)$ is the complex, frequency-dependent,
viscoacoustic modulus.

The dissipative aspect of a material is often described by the
quality factor $Q$, defined as the ratio of the real and imaginary
parts of the viscoacoustic modulus. It expresses how attenuating a
material is and corresponds to the number of wavelengths a wave
can propagate through the medium before its amplitude has
decreased by $e^{-
  \pi}$,
\begin{equation}
\dsp Q(\omega)=\frac{\Re(\mu(\omega))}{\Im(\mu(\omega))}=\frac{1}{\tan(\phi(\omega))},
\label{definitionQ}
\end{equation}
where $\phi(\omega)$ is the phase of $\mu(\omega)$.

In seismic applications, $Q$ is usually assumed to be frequency-
independent or only slowly frequency-dependent. In this case (i.e.
when $Q$ is constant in frequency), the viscoacoustic modulus is
given analytically by  Kjartansson's model \cite{Kjar1979},
\begin{equation}
\dsp \mu(\omega) = \mu_{ref} \left( \frac{\mathbf{i} \omega}{\omega_{ref}}
\right) ^{\frac{2}{\pi} arctan(Q^{-1})}.
\label{Kjar}
\end{equation}
This analytical formulation will be useful for validation of the
numerical results in the next sections.

\begin{figure}[htbp]
\centering\psfig{figure=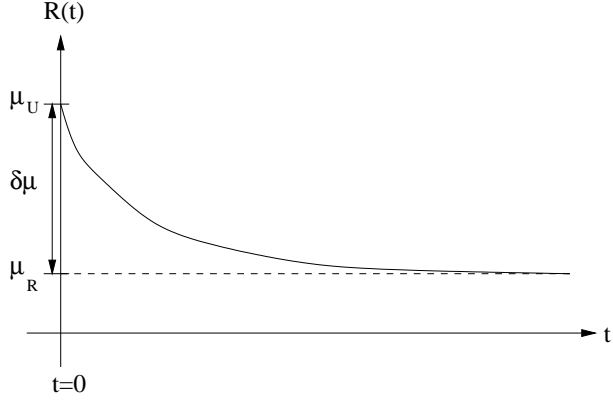,width=8cm}
\caption{Schematic example of the relaxation function $R(t)$}
\label{fig:relax}
\end{figure}

In the time domain, the constitutive relation
(\ref{relationpufreq}) is expressed in terms of a convolution
operator, denoted here by $\star_{t}$,
\begin{equation}
\dsp p(t)=\mu(t) \star_{t}  \Div \mathbf{u}(t).
\label{relationputime}
\end{equation}
The discretization of this equation requires saving in memory the
whole history of the solution at all points of the computational
domain and is thus very expensive. To overcome this inconvenience,
we approximate the viscoacoustic modulus by a rational function in
frequency, as was proposed in \cite{DM84,
  EK87,CKK88a, BRS95}. It is convenient in the following to
introduce the relaxation function $R( \mathbf{x},t)$, defined by, (see
Figure \ref{fig:relax}),
\begin{equation}
\dsp \mu(\mathbf{x},t)=\frac{\partial R(\mathbf{x},t)}{ \partial t} \text{
  ; }R(\mathbf{x},t)=\left(\mu_{R} (\mathbf{x})+ \delta
\mu(\mathbf{x}) \int_{0}^{+\infty}r(\mathbf{x},\omega') e^{-\omega't}
d\omega'\right)H(t),
\label{Relaxfun}
\end{equation}
where $\mu_{R}$ is the relaxed modulus,
$$
\dsp \mu_{R}(\mathbf{x})=\lim_{t\rightarrow +\infty}R(t),
$$
$\mu_{U}$ is the unrelaxed modulus,

$$\dsp \mu_{U}(\mathbf{x})=\mu_{R}(\mathbf{x})+ \delta
\mu(\mathbf{x})  =\lim_{t\rightarrow 0}R(t),
$$
$r(\mathbf{x},\omega')$ is the normalized relaxation spectrum
satisfying
$$
\dsp \int_0^{+\infty}r(\mathbf{x},\omega') d\omega'=1,
$$
and $H(t)$ is the Heaviside function.

Using $\mu(\mathbf{x},t)$ defined by (\ref{Relaxfun}) in
(\ref{relationputime}) gives, $$ \dsp p(\mathbf{x},t)=\mu_{U}
(\mathbf{x}) div \mathbf{u}(\mathbf{x},t) - \delta \mu
(\mathbf{x}) \int_{-\infty}^{t}  \int_{0}^{+\infty}
\omega'r(\mathbf{x},\omega') e^{-\omega'(t-\tau)}div
\mathbf{u}(\mathbf{x},\tau) d\omega' d \tau .$$

We now assume that the relaxation spectrum can be discretized by
$L$ single peaks of amplitude $\alpha_{l}$ at relaxation
frequencies $\omega_{l}$, $l \in [1...L]$,

$$\dsp r(\mathbf{x},\omega)=\sum_{l=1}^{L} \alpha_{l}(\mathbf{x}) \delta
(\omega-\omega_{l}(\mathbf{x}))\text{ ; }\sum_{l=1}^{L}
\alpha_{l}(\mathbf{x})=1,
$$
In this case, we get,
$$
\dsp R(\mathbf{x},t) \simeq R_{l}(\mathbf{x},t)=\left(\mu_{R}
(\mathbf{x})+\sum_{l=1}^{L} \alpha_{l}(\mathbf{x})
e^{-\omega_{l}(\mathbf{x})t} \right)H(t) ,
$$
and
\begin{equation}
\dsp \mu(\mathbf{x},\omega)\simeq \mu_{l}(\mathbf{x},\omega)=\mu_{R}
(\mathbf{x})\left(1+\sum_{l=1}^{L}
\frac{y_{l}(\mathbf{x})\mathbf{i}\omega}{\mathbf{i}\omega+\omega_{l}(\mathbf{x})}
\right).
\label{Modviscoap}
\end{equation}
In (\ref{Modviscoap}), we introduced $y_{l}(\mathbf{x})$ defined by,

$$
\dsp y_{l}(\mathbf{x})=\frac{\delta
  \mu(\mathbf{x})}{\mu_{R}(\mathbf{x})}\alpha_{l}(\mathbf{x}) \text{,
  with the normalization relation}\sum_{l=1}^{L}y_{l}(\mathbf{x})=\frac{\delta
  \mu(\mathbf{x})}{\mu_{R}(\mathbf{x})}.
$$ Notice that equation (\ref{Modviscoap}) can be obtained if one
assumes that $\mu(\mathbf{x},\omega)$ can be approximated by a
rational function of $(\mathbf{i}\omega)$,
\begin{equation}
\dsp \mu(\mathbf{x},\omega) \simeq
\mu_{l}(\mathbf{x},\omega)=\frac{P_{L}(\mathbf{x},\mathbf{i}\omega)}{Q_{L}(\mathbf{x},\mathbf{i}\omega)},
\label{ModPol}
\end{equation}
with $P_{L}$ and $Q_{L}$ being polynomials of degree $L$ in
$(i\omega)$. Then (\ref{Modviscoap}) can be interpreted as an expansion
of (\ref{ModPol}) into partial fractions \cite{EK87}. Thus
approximating the viscoacoustic modulus by a rational function is
equivalent to approximating the relaxation spectrum by a discrete one.

For computational reasons, it is natural to search for rational
function approximations of the viscoacoustic modulus, which
minimize the ratio: number of unknowns/accuracy. We therefore
address in the following the question of finding an accurate
low-order approximation of the viscoacoustic modulus.
\section{Approximation of the viscoacoustic modulus}\label{s2}
We now briefly introduce the different
approximation methods previously proposed in the literature.
\subsection{Pad\'e approximation method}{\label{s21}}
The use of the simple Pad\'e approximation in the framework of
viscoacoustic wave propagation was proposed in \cite{DM84}. Letting
$\dsp z=-\frac{1}{\mathbf{i}\omega}$ and introducing,
$$\dsp
\chi(\mathbf{x},z)=\int_0^{+\infty}\frac{\omega'r(\mathbf{x},\omega')}{1-\omega'z}d\omega',$$
$\mu(\mathbf{x},\omega)$ can be re-written in the following form,
$$
\dsp \mu(\mathbf{x},\omega)=\mu_{U}(\mathbf{x})+\delta \mu(\mathbf{x})
 z \chi(\mathbf{x},z).
$$ The Pad\'e approximation is then used for expanding
$\chi(\mathbf{x},z)$ into a rational function with numerator of
degree $L-1$ and denominator of degree $L$. Using the well-known
(\cite{Sze1939, Bre1980}) relations between Pad\'e approximations
and orthogonal polynomials one gets,

$$\dsp
\chi(\mathbf{x},z)=\sum_{l=1}^{L}\frac{\lambda_{l}(\mathbf{x})}{1-\omega_{l}(\mathbf{x})z},
$$ where $\omega_{l}$ are the zeros of the orthogonal polynomial
$P_{L}$, and $\lambda_{l}$ are the residuals given by, $$ \dsp
\lambda_{l}(\mathbf{x})=\frac{k_{L}}{k_{L-1}P_{L-1}(\omega_{l}(\mathbf{x}))P'_n(\omega_{l}(\mathbf{x}))},
$$ $k_{L}$ being the leading coefficient of $P_{L}$ and where the
prime denotes the derivative of $P_{L}$. Recall that the
orthogonal polynomials are defined by,

$$\dsp
\int_{\Omega1}^{\Omega2}P_{n}(\omega')P_{m}(\omega')\omega'r(\omega')d\omega'=\delta_{mn},
$$

where $\delta_{mn}$ is the Kronecker symbol. When the quality
factor is constant over a frequency band, $\lambda_{l}$ and
$\omega_{l}$ can be obtained in closed form. Moreover, when $Q \ll
1$, the relaxation spectrum $r(\mathbf{x},\omega)$ is proportional
to $\omega^{-1}$. Assuming that $r(\mathbf{x},\omega)$ is zero
outside the frequency interval $[\Omega_{1},\Omega_{2}]$ we obtain
the approximation,
\begin{equation}
\dsp \mu(\mathbf{x},\omega) \simeq \mu_{l}(\omega)=\mu_{U} \left(
1-\frac{\Omega2-\Omega1}{\pi
  Q}\sum_{l=1}^{L}\frac{\nu_{l}(\mathbf{x})}{i\omega+\omega_{l}(\mathbf{x})} \right),
\label{ApPade}
\end{equation}
where $\dsp
\omega_{l}=\frac{1}{2}[x_{l}(\Omega2-\Omega1)+\Omega2+\Omega1]$,
$x_{l}$ and $\nu_{l}$ being respectively the zeros and weights of
the Legendre polynomials. Notice that the relaxation frequencies
$\omega_{l}$ are in this case equidistant on a linear scale. The
main advantage of this approximation is that all data are
analytically determined. For more details on this method the
reader can refer to \cite{DM84}.
\subsection{Generalized Maxwell Body approximation method}{\label{s22}}
We describe here the method proposed in \cite{EK87}. First let us re-write
(\ref{Modviscoap}) as,
\begin{equation}
\dsp \mu_{l}(\mathbf{x},\omega)=\mu_{R}(\mathbf{x})+\delta
\mu(\mathbf{x})\sum_{l=1}^{L}\frac{\alpha_{l}(\mathbf{x}) \mathbf{i}
  \omega}{\mathbf{i}\omega+\omega_{l}(\mathbf{x})}.
\label{ModPo2}
\end{equation}
Each term of (\ref{ModPo2}) can be interpreted as a classical Maxwell
body with viscosity $\dsp \alpha_{l} \frac{\delta \mu}{\omega_{l}}$ and
elastic modulus $\alpha_{l} \delta \mu$. The term $\mu_{R}$ in (\ref{ModPo2})
represents an additional elastic element. The $Q$-law for the generalized
Maxwell body approximation can be obtained from (\ref{ModPo2}),
\begin{equation}
\dsp Q(\mathbf{x},\omega)^{-1}=\frac{\Im(\mu(\mathbf{x},\omega))}
     {\Re(\mu(\mathbf{x},\omega))}=\frac{\sum_{l=1}^{L}y_{l}(\mathbf{x})\frac{\frac{\omega}{\omega_{l}(\mathbf{x})}}
     {1+(\frac{\omega}{\omega_{l}(\mathbf{x})})^2}}
     {1+\sum_{l=1}^{L}y_{l}(\mathbf{x})\frac{(\frac{\omega}{\omega_{l}
(\mathbf{x})})^2}{1+(\frac{\omega}{\omega_{l}(\mathbf{x})})^2}}.
\label{QMax}
\end{equation}
Assuming now that $\delta \mu \ll \mu_{R}$, (\ref{QMax}) becomes,
\begin{equation}
\dsp Q(\mathbf{x},\omega)^{-1} \backsimeq \frac{\delta
  \mu(\mathbf{x})}{\mu_{R}(\mathbf{x})}
  \sum_{l=1}^{L}\alpha_{l}(\mathbf{x})
  \frac{\frac{\omega}{\omega_{l}(\mathbf{x})}} {1+(\frac{\omega}
  {\omega_{l}(\mathbf{x})})^2}.
\label{QMax2}
\end{equation}
This means that $Q(\omega)^{-1}$ is approximately the sum of $n$
Debye functions with maxima $\dsp \alpha_{l} \frac{\delta
\mu}{2\mu_{R}}$ located at frequencies $\omega_{l}$. If $Q$ is
fairly constant in a frequency band, the most natural choice for
the relaxation frequencies $\omega_{l}$ is a logarithmic
equidistant distribution. In this case, to obtain a good
approximation of $Q(\omega)^{-1}$, the distance between two
adjacent relaxation frequencies should be chosen smaller or equal
to the half-width of the Debye function (1.144 decades). In
\cite{EK87} two ways for choosing $\omega_l$ were proposed:
$\omega_{l}$ can be chosen logarithmically-equidistant in the
frequency band $[\Omega1,\Omega2]$ or determined by $\dsp
\omega_{l}=\frac{2\omega_{dom}}{10^{l}}$ where $\omega_{dom}$ is
the dominant (central) frequency of the source considered in the
simulations. In both cases, the coefficients $y_{l}$ are obtained
by solving the overdetermined linear system
\begin{equation}
\dsp \sum_{l=1}^{L}y_{l}(\mathbf{x})\widetilde{\omega}_{k}
(\mathbf{x})\frac{\omega_{l}
  (\mathbf{x})-\widetilde{Q}^{-1}(\mathbf{x}, \widetilde{\omega}_{k}
  (\mathbf{x}))\widetilde{\omega}_{k} (\mathbf{x})}{\omega_{l}
  (\mathbf{x})^2+\widetilde{\omega}_{k}
  (\mathbf{x})^2}=\widetilde{Q}^{-1}
(\mathbf{x},\widetilde{\omega}_{k} (\mathbf{x}))\textit{, k }\in[1,2,..K],
\label{Dety_l}
\end{equation}
where, $\widetilde{\omega}_{k}$ are defined by
$$
\barr{l}
\widetilde{\omega}_{1}=\Omega_{1},\\[8pt]
\dsp \widetilde{\omega}_{k+1}=\widetilde{\omega}_{k} (
\frac{\Omega2}{\Omega1} ) ^\frac{1}{2}.
\earr
$$
Let us remark that the determination of $\omega_{l}$ for this
approximation is based on an empirical study.
\subsection{The $\tau$-method}{\label{s23}}
This method, proposed in \cite{BRS95} is based on the observation
that dissipation due  to only one ``Maxwell Body'' can be
determined by a unique dimensionless parameter $\tau$. More
precisely, for $Q \gg 1$  and $L=1$, equation (\ref{QMax2})
becomes, $$ \dsp
Q(\mathbf{x},\omega)^{-1}=\frac{\frac{\omega}{\omega_{1}(\mathbf{x})}
  \tau (\mathbf{x})}{1+(\frac{\omega}{\omega_{1}(\mathbf{x})})^{2}},
$$

where $\tau = y_1 \ll 1$. It is then easy to see (cf. \cite{BRS95}),
that $\omega_{1}$ essentially determines the
frequency behavior of $Q$ while $\tau$ determines its magnitude. In
the general case for $L>1$, and when one seeks an approximation of a
constant $Q$ value, $y_{l}$ are quasi-constant

and equation (\ref{QMax2}) can be approximated by,
\begin{equation}
\dsp Q(\mathbf{x},\omega)^{-1}=
\sum_{l=1}^{L}\frac{\frac{\omega}{\omega_{l}(\mathbf{x})}
  \tau(\mathbf{x})}{1+(\frac{\omega}{\omega_{l}(\mathbf{x})})^{2}}.
\label{tauL}
\end{equation}
In (\ref{tauL}), $Q(\omega)^{-1}$ is linear in $\tau$. One can
therefore find the best approximation, in the least-squares sense,
over a predefined frequency range to any $Q_{0}$ by minimizing over
$\tau$ the expression,
\begin{equation}
J=\int_{\Omega_1}^{\Omega_2}(Q^{-1}(\omega,\omega_{l},\tau)-Q_0^{-1})^2
d\omega.
\label{Mintau}
\end{equation}
The approximation of the viscoacoustic modulus in this case is,
\begin{equation}
\mu_{l}(\mathbf{x},\omega)=\mu_{R} (\mathbf{x})\left(1+\sum_{l=1}^{L}
\frac{\tau(\mathbf{x})\mathbf{i}\omega}{\mathbf{i}\omega+\omega_{l}(\mathbf{x})} \right).
\label{Modviscotau}
\end{equation}
The relaxation frequencies $\omega_{l}$ are chosen, as for the
``Generalized Maxwell Body'' method, equidistant on a logarithmic
scale. Equation (\ref{Modviscotau}) leads in general to an over-estimation of
the value of $Q$.  Thus the authors in \cite{BRS95} suggest to use in
the definition of $J$ (\ref{Mintau}) a
value for $Q_0$ slightly smaller than the desired one.
 This value is also  chosen empirically.

\subsection{Comparison of the different approximation methods}{\label{s24}}
To test the accuracy of the different approximation methods
previously presented, we compute the response of a one-dimensional
viscoacoustic homogeneous medium to the following pulse,
\begin{equation}
\dsp s(t) = {\sin \left( \frac{2\pi t}{T} \right) }-0.5{\sin \left(\frac{4\pi t}{T} \right)} \text{  for  }
0<t<T \text{,  } T=0.3s.
\label{Input}
\end{equation}
The solution is obtained by convolving the source function $s(t)$ with
the dissipation operator $D(t)$ (the Green's function for the 1D
problem). For an
arbitrary dissipation law, the Fourier transform $D(\omega)$ of $D(t)$
is given by \cite{EK87},
\begin{equation}
\dsp D(\omega)=e^{\mathbf{i}\omega t^{\star} Q(\omega_{r})
  \left(1-\frac{c(\omega_{r})}{\nu(\omega)}\right)},
\label{OpDw}
\end{equation}
where $c(\omega_{r})$ is the phase velocity at the reference frequency
$\omega_r$, $\nu(\omega)$ the complex velocity, and
$\dsp t^\star=\frac{x}{c(\omega_r)Q(\omega_r)}$ the dissipation
time.  For a frequency independent $Q$,
the value of $\dsp
\frac{c(\omega_{r})}{\nu(\omega)}=\frac{|\mu(\omega_{r})|}{\mu(\omega)}$,

can be determined from equation (\ref{Kjar}) combined with one of
(\ref{Modviscoap}), (\ref{ApPade}) or (\ref{Modviscotau}), depending
on the approximation method used.

In the numerical example, we want to approximate $Q=20$ over the
frequency range $[10^{-2}, 10^{2}]$ Hz, like in \cite{EK87}.

\begin{figure}[htbp]
\centering\psfig{figure=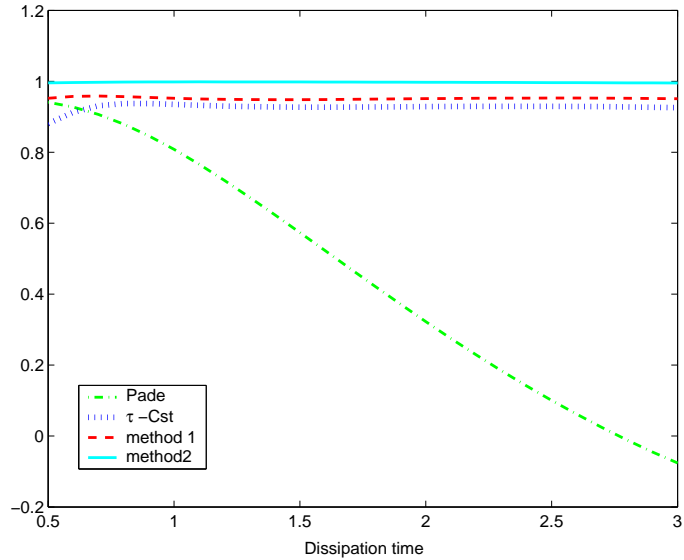,width=9cm}
\caption{Comparison between the different approximation
  methods. Correlation coefficient between the exact and the
  approximated solutions as a function of the dissipation time.}
\label{Cor}
\end{figure}

To better illustrate the results, we present in Figure \ref{Cor} the
evolution of the correlation coefficient between the exact solution
(the one obtained for the viscosity modulus calculated from
(\ref{Kjar})) and the different approximated ones (calculated with the
viscosity modulus provided by (\ref{Modviscoap}),
(\ref{ApPade}) or (\ref{Modviscotau})) as a function of
the dissipation time.

More precisely in Figure \ref{Cor} we compare the results obtained with the
following approximations,
\begin{itemize}
\item Pad\'e approximation with $L=5$.
\item Maxwell Body approximation with $L=3$ and relaxation frequencies
  chosen logarithmically-equidistant over the frequency range
  $[10^{-1.5},10^{1.5}]$ Hz (cf. \cite{EK87}). We call this choice
  method~1.
\item Maxwell Body approximation with $L=3$ and relaxation frequencies
  chosen equidistant on a logarithmic scale, such that,
  $\omega_{l}=\frac{2\omega_{dom}}{10^{l}}$ (cf. \cite{EK87}). We call
  this choice  method~2.
\item The $\tau$-method with $Q_{0}=17.6$ (value proposed in
  \cite{BRS95} to model the propagation in a viscoacoustic medium with
  $Q=20$) and $L=3$.
\end{itemize}
The results illustrated in Figure \ref{Cor}, show that the Pad\'e
approximant provides good accuracy only for short dissipation
times, as demonstarted in \cite{EK87}. The $\tau$-method provides
a good  accuracy/number of calculations ratio. However, we did not
choose this method because $Q_{0}$ has to be calibrated
empirically in order to get good results. The ``Generalized
Maxwell Body'' approximation method seems to be a good compromise
between accuracy, number of calculations, and implementation
simplicity. As our aim is to simulate viscoacoustic wave
propagation in heterogeneous media for large dissipation times, we
chose a  method which is a hybrid of the Maxwell approximation
methods 1 and 2.
\subsection{Proposed method}{\label{s25}
In practice, the source type used depends on the application of
interest. In our case, the main application of interest is
seismic wave propagation for which a Ricker wavelet is often used as
source function,

\begin{equation}
\label{ricker}
f(t) = -2 \alpha^{2} \left(1-2\alpha^{2}\left(t-\frac{1}{f_{0}}
\right)^{2} \right){\exp\left(-\alpha^{2}\left(t-\frac{1}{f_{0}}
  \right) \right)} \text{,  for  }0<t\leq\frac{2}{f_{0}}\text{,  with
}\alpha=\pi f_{0}.
\end{equation}
\begin{figure}[htbp]
\begin{minipage}{7.5cm}
\centering\psfig{figure=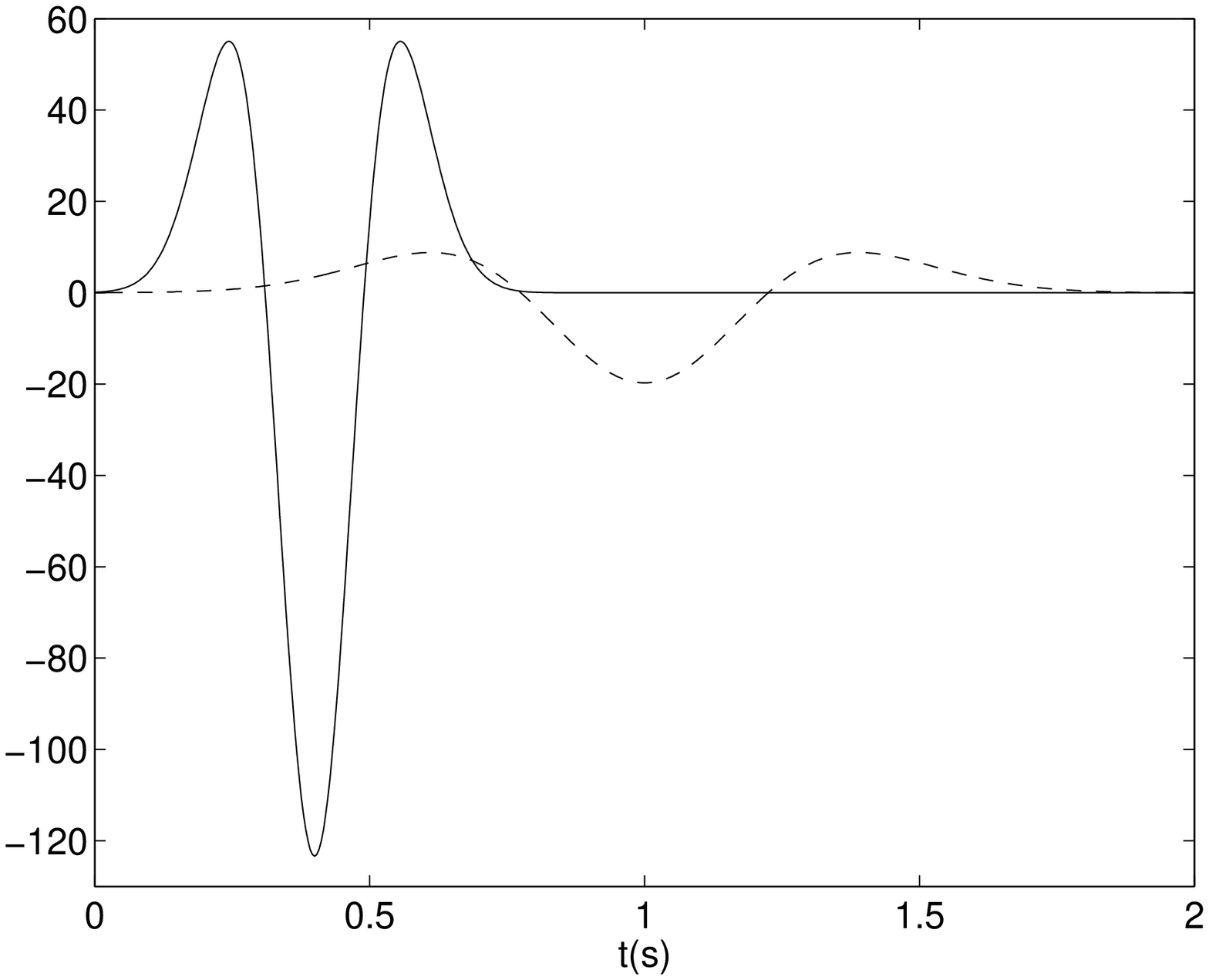,width=8cm}
\caption{The Ricker wavelet $f(t)$ for $f_{0}=2.5Hz$ (solid line) and for
  $f_{0}=1Hz$ (dashed line).}
\label{source1}
\end{minipage} \hfill
\begin{minipage}{7.5cm}
\centering\psfig{figure=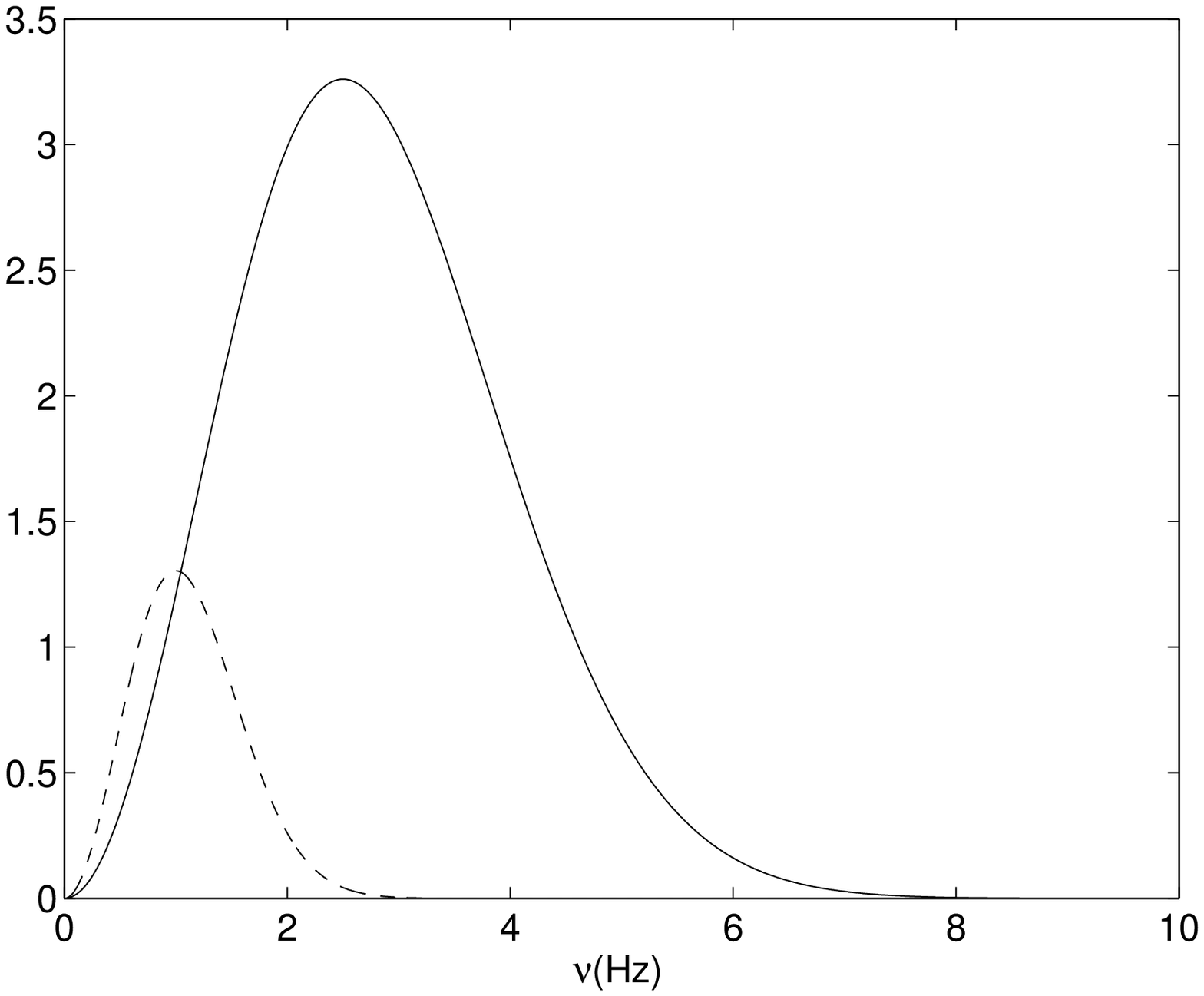,width=8cm}
\caption{The frequency spectrum of the Ricker wavelet for
  $f_{0}=2.5Hz$ (solid line) and for $f_{0}=1Hz$ (dasged line).}
\label{source2}
\end{minipage}
\end{figure}
In Figures \ref{source1} and \ref{source2} we display the source
function (\ref{ricker}) and its spectrum for two different values
of $f_0$. Compared to the source function defined by
(\ref{Input}), the Ricker wavelet has a broader frequency spectrum
and method 2 did not give as good results in this case as the ones
obtained with the source (\ref{Input}). Following the ideas in
\cite{EK87}, we want to find a way to choose the frequency band
$[\Omega1, \Omega2]$ as a function of the source type and then
determine the relaxation frequencies $\omega_{l}$ logarithmically
equidistant in this bandwidth. We found that a good choice for a
Ricker wavelet type of source is $\dsp [\Omega1,
\Omega2]=[\frac{\omega_{max}}{100},\omega_{max}]$, where
$\omega_{max}$ is the maximal frequency of the employed source
spectrum. None of the above approximation methods is completely
satisfactory in our opinion because the choice of the relaxation
frequencies is always empirical.

To avoid this, one can follow the approach proposed in
\cite{Asv2003} where a non-linear minimization problem is
considered which permits to determine all the coefficients (both
$y_{l}$ and $\omega_{l}$ $\forall l \in [1, L]$). However, this
method is more expensive and although it improves the accuracy of
the solution for media with high damping ($Q\le 10$) it  provides
quite similar results with the proposed method for propagation in
weakly attenuating media ($Q \ge 10$) \cite{Asv2003}. As we are
interested in media with quality factors greater than $10$, we
will use in the following the linear minimization method (system
(\ref{Dety_l})).

\section{The mixed velocity-pressure formulation}{\label{s3}}
By incorporating (\ref{Modviscoap}) into (\ref{relationpufreq})
we get,
\begin{equation}
\dsp
p(\mathbf{x},\omega)=
\mu_{R}(\mathbf{x})\Div(\mathbf{u}(\mathbf{x},\omega))+\mu_{R}(\mathbf{x})
\sum_{l=1}^{L}\frac{y_{l}(\mathbf{x})\mathbf{i}\omega}{\mathbf{i}\omega
  +\omega_{l}
  (\mathbf{x})} \Div(\mathbf{u}(\mathbf{x},\omega)).
\label{Pap2}
\end{equation}
We now introduce the memory variables $\eta_{l}$ defined by,
\begin{equation}
\label{eta1}
\dsp (\mathbf{i}\omega+\omega_{l}(\mathbf{x}))\eta_{l}(\mathbf{x},\omega) =
  \mu_{R}(\mathbf{x})y_{l}(\mathbf{x}) \Div(\mathbf{v}(\mathbf{x},\omega)),
\end{equation}
where $\mathbf{v}$ is the velocity, i.e., the time derivative of
the displacement $\mathbf{u}$. Equation (\ref{eta1}) in the time
domain becomes,
\begin{equation}
\dsp \frac{\partial \eta_{l}(\mathbf{x},t)}{\partial t} +
\omega_{l}(\mathbf{x})\eta_{l}(\mathbf{x},t) =
\mu_{R}(\mathbf{x})y_{l}(\mathbf{x}) \Div(\mathbf{v}(\mathbf{x},t)).
\label{etafin}
\end{equation}
Using the definition of $\eta_{l}$ and multiplying (\ref{Pap2}) by
$(\mathbf{i}\omega)$, we get,
$$
\dsp (\mathbf{i}\omega)p(\mathbf{x},\omega) =
  \mu_{R}(\mathbf{x})\Div(\mathbf{v}(\mathbf{x},\omega)) +
  \sum_{l=1}^{L}(\mathbf{i}\omega)\eta_{l}(\mathbf{x},\omega),
$$ or equivalently in the time domain,
\begin{equation}
\displaystyle \frac{\partial p}{\partial t} = \mu_{R} \Div(\mathbf{v})
+ \sum_{l=1}^{n}
\frac{\partial \eta_{l}}{\partial t}.
\label{pfin}
\end{equation}
Combining (\ref{pfin}), (\ref{etafin}) and the equation of motion,
we obtain our final system of equations,
\begin{equation}
\left\{ \begin{array}{ll}
\displaystyle \rho \frac{\partial \mathbf{v}}{\partial t}-\nabla p =
\mathbf{f}&\text{in }\Omega \times ]0,T], \\[12pt]
\displaystyle \frac{\partial p}{\partial t} - \sum_{l=1}^{n}
\frac{\partial \eta_{l}}{\partial t} = \mu_{R} \Div(\mathbf{v}) &
\text{in } \Omega \times ]0,T], \\[12pt]
\displaystyle \frac{\partial \eta_{l}}{\partial t} +
\omega_{l}\eta_{l} = \mu_{R}y_{l} \Div(\mathbf{v}), \forall l & \text{in
} \Omega \times ]0,T].
\label{systime}
\end{array}
\right.
\end{equation}
Equivalently, one can chose to eliminate the pressure and obtain a
second-order-in-time equation for the displacement by introducing
adequate memory variables \cite{EK87}. We prefer, however, the
first-order velocity-pressure formulation for the following
reasons,
\begin{itemize}
\item It can be coupled with the fictitious domain method for taking
  into account diffraction by objects of complicated geometry.
\item A perfectly matched layer model (PML) can be written for this
  system. This permits us to simulate efficiently wave propagation in
  unbounded domains.
\item This system is easier to implement in heterogeneous media, since
  it does not require an approximation of the spatial derivatives of the
  physical parameters.
\end{itemize}
An equivalent first-order velocity-pressure system is proposed in
\cite{CKK88a} and \cite{BRS95}. In \cite{CKK88a} the authors used
a pseudospectral method for the discretization while in
\cite{BRS95} a staggered finite difference scheme was used. Our
aim being to couple this system with the fictitious domain method,
we propose here instead the use of a mixed-finite element method
on regular grids. A similar approach was proposed in \cite{BEJ03}
where the authors use a mixed-finite element method to discretize
the velocity-stress formulation for viscoelastic wave propagation.
\section{Discretisation}{\label{s4}}
A mixed formulation associated to equations (\ref{systime}) is given by,
\begin{equation}
\left\{
\begin{array}{ll}
\text{Find } (\mathbf{v},p,H) :  ]0,T[  \longmapsto X   \times M
    \times (M)^{L} \text{ s.t. :} & \\[6pt]
\displaystyle \frac{d}{dt} (\rho \mathbf{v},\mathbf{w}) +
    b(\mathbf{w},p)  = (\mathbf{f},\mathbf{w}) , & \displaystyle
    \forall {\mathbf w} \in X, \\ [12pt]
\displaystyle \frac{d}{dt} (\frac{1}{\mu_{R}} p,q) - \sum_{l=1}^{L}
    \frac{d}{dt}(\frac{1}{\mu_{R}} \eta_{l},q) - b({\mathbf{v}},q)=0,
    & \displaystyle \forall q \in M , \\[12pt]
\displaystyle \frac{d}{dt} (\frac{1}{\mu_{R} y_{l}} \eta_l,q)
+ (\frac{\omega_{l}}{\mu_{R} y_{l}} \eta_{l},q) -
    b({\mathbf{v}},q)=0,\ \forall l,
& \displaystyle \forall q \in M ,
\end{array}
\right.
\label{Varfor}
\end{equation}
where $H$ is the L-dimensional vector with components $\eta_l$, and

$$\begin{array}{ll}
\displaystyle b(\mathbf{w},q) = \int_{\Omega}  q \ \Div\mathbf{w}\
d\mathbf{x},
& \displaystyle  \forall (\mathbf{w},q) \in
X \times M .
\end{array}
$$ The functional spaces are $X = H(\Div;\Omega)$, and $M =
L^2(\Omega)$.\\ We now introduce  some finite element spaces
$X_{h} \subset X$, and $M_{h} \subset M$ of dimensions $N_1$ and
$N_2$ respectively. The semi-discretization in the space of
problem (\ref{Varfor})  is,
\begin{equation}
\left\{
\begin{array}{l}
(V_{h},P_{h},H_{h}) \in L^{2}(0,T;I\!R^{N_{1}}) \times
L^2(0,T;I\!R^{N_{2}}) \times L^{2}(0,T;(I\!R^{N_2})^{L})
\text{ s.t. :}
\\[6pt]
\displaystyle M_{v} \frac{d V_{h}}{dt} + B_{h} P_{h} = F_{h} ,
\\[6pt]
\displaystyle M_{p} \frac{d P_{h}}{dt} - \sum_{l=1}^{L}  M_{p} \frac{d
  (H_{h})_{l}}{dt} -  B_h^T V_{h} = 0 ,
\\[6pt]
\displaystyle M_{y} \frac{d (H_{h})_{l}}{dt} +  M_{\omega} (H_{h})_{l}
-B_h^T V_{h} = 0 ,\ \forall l,
\end{array}
\right.
\label{SemDi}
\end{equation}
where $B_{h}^{T}$ denotes the transpose of $B_{h}$.\\ In practice,
we only consider regular domains in $I\!R^{d}$, $d=1,2$ that can
be discretized with a uniform mesh ${\cal T}_{h}$ composed by
segments or squares of size $h$, depending on the dimension of the
problem. The finite element spaces we use are $$ \displaystyle
X_{h} = \left\{ \mathbf{w}_{h} \text{ in } X
  \left/ \right. \forall  K
\in  {\cal T}_{h},  \mathbf{w}_{h}\left|_K \right.
\in (Q_1)^d \right\},$$
and $$M_h = \left\{ {q}_h \in L^2 \left/ \right. \forall  K
\in  {\cal T}_{h},  {q}_{h\left|_{K}
    \right.} \in {P_0}(K) \right\}. $$
This mixed finite element was introduced in \cite{BJT1-SIAM} and is
illustrated in Figure \ref{VP0}.\\
\begin{figure}[htbp]
\centering\psfig{figure=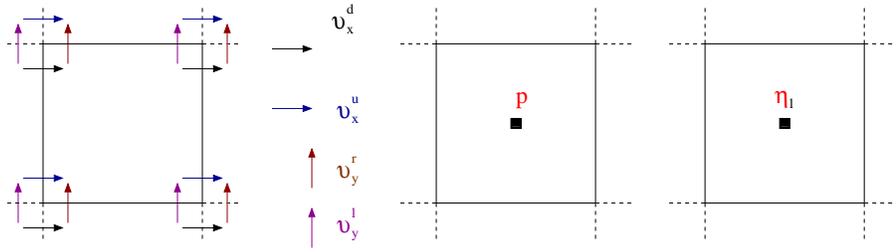,width=12cm}
\caption{Finite element $v_{h} \in  X_{h}, (p_{h}, (\eta_{l})_{h}) \in M_{h} \times M_{h}$}
\label{VP0}
\end{figure}
When coupled with the fictitious domain method, this choice of
finite elements presents some inconveniences. In particular, for
the acoustic wave equation problem we cannot prove the convergence
of the method from the theoretical point of view. Moreover,
numerical results show that the method converges under restrictive
conditions on the discretization parameters. Thus, when the method
is coupled with the fictitious domains, we replace $M_h$ by
$M_h^1$ defined by, $$M_h^1 = \left\{ {q}_h \in L^2   \left/
\right. \forall  K \in  {\cal T}_{h},  {q}_{h\left|_{K}
    \right.} \in {P_1}(K) \right\}.$$
This finite element is presented Figure \ref{Vp1}. In this case
convergence for the acoustic waves problem coupled with
the fictitious domain method was obtained \cite{BRT}.

For computational reasons, however, it is natural to seek a
discretization which uses the least number of variables. In the
proposed scheme pressure is thus discretized on the space
$M_{h}^{1}$ and the memory  variables $\eta_{l}$ are discretized
on $M_{h}$. The semi-discretization of the problem (\ref{Varfor})
in this case is,
\begin{equation}
\left\{
\begin{array}{l}
(V_{h},P_{h},H_{h}) \in L^{2}(0,T;I\!R^{N_{1}}) \times
L^2(0,T;I\!R^{N_{2}}) \times L^{2}(0,T;(I\!R^{N_3})^{L})
\text{ s.t. :}
\\[6pt]
\displaystyle M_{v}^{1} \frac{d V_{h}}{dt} + B_{h}^{1} P_{h} = F_{h}^{1} ,
\\[6pt]
\displaystyle M_{p}^{1} \frac{d P_{h}}{dt} - \sum_{l=1}^{L}  M_{p}^{0} \frac{d
  (H_{h})_{l}}{dt} -   B_{h}^{1,T} V_{h} = 0 ,
\\[6pt]
\displaystyle M_{y}^{0} \frac{d (H_{h})_{l}}{dt} +  M_{\omega}^{0} (H_{h})_{l}
-  B_h^{0,T} V_{h} = 0 ,\ \forall l.
\end{array}
\right.
\label{SemDi1}
\end{equation}
\begin{figure}[htbp]
\centering\psfig{figure=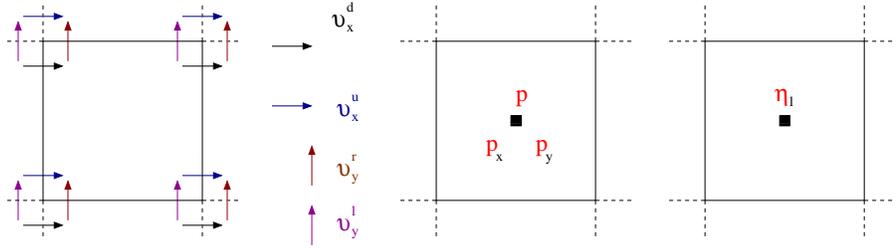,width=12cm}
\caption{Finite element $v_{h} \in  X_{h},  (p_{h}, (\eta_{l})_{h}) \in M_{h}^{1} \times M_{h} $}
\label{Vp1}
\end{figure}
In both cases (pressure discretized on $M_{h}$ or $M_{h}^{1}$), we
use a second order centered finite difference scheme for the
discretization in time (here presented in the more general case
with the pressure discretized on $M_{h}^{1}$),
\begin{equation}
\left\{
\begin{array}{l}
(V_{h}^{n+1}, P_{h}^{n+\frac{3}{2}}, H_{h}^{n+\frac{3}{2}})
\in I\!R^{N_1} \times I\!R^{N_2}\times (I\!R^{N_3})^{L},
 \\[6pt]
\displaystyle  M_{v}^{1}\frac{V_{h}^{n+1}-V_{h}^{n} }{\Delta t}  +
B_{h}^{1} P_{h}^{n+\frac{1}{2}}
= \left(F_{h}^{1}\right)^{n+1/2}, \\[6pt]
\displaystyle  M_{p}^{1} \frac{P_{h}^{n+\frac{3}{2}}-P_{h}^{n+\frac{1}{2}}
}{\Delta t}  -
\sum_{l=1}^L  M_{p}^{0}
\frac{(H_{h})_{l}^{n+\frac{3}{2}}-(H_{h})_{l}^{n+\frac{1}{2}} }{\Delta
  t}
-  B_{h}^{1,T}  V_{h}^{n+1} =0, \\[6pt]
\displaystyle M_{y}^{0}
\frac{(H_{h})_{l}^{n+\frac{3}{2}}-(H_{h})_{l}^{n+\frac{1}{2}} }{\Delta
  t}
+ M_{\omega}^{0} \frac{(H_{h})_{l}^{n+\frac{3}{2}}+(H_{h})_{l}^{n+\frac{1}{2}} }{2}
- B_{h}^{0,T} V_{h}^{n+1} = 0 , \forall l.
\end{array}
\right.
\label{TotDi}
\end{equation}
The numerical scheme (\ref{TotDi}) becomes explicit in time when
an adequate quadrature formula is used to approximate the matrix
$M_{v}^1$. Note that the other mass matrices ($M_{p}^0$,
$M_{p}^1$, $M_{y}^0$, and $M_{\omega}^0$) are diagonal, since the
spaces $M_{h}$ and $M_{h}^1$ are composed of discontinuous
functions. For more details on the quadrature formulas used we
refer the reader to \cite{BJT1-SIAM}.
%
\section{Stability and dispersion analysis}{\label{s5}}
For the continuous problem, the energy is defined by
\begin{equation}
\varepsilon= \frac{1}{2} \left(\rho
\mathbf{v},\mathbf{v}\right)+\frac{1}{2}\left(p-\sum_{l=1}^{L}\eta_{l},
p-\sum_{l=1}^{L}\eta_{l} \right)+ \sum_{l=1}^{L}\frac{1}{2y_{l}
  \mu_{R}} \left(\eta_{l}, \eta_{l} \right).
\label{Encont}
\end{equation}
This quantity is positive (for $y_l$ positive) and we have,
\begin{equation}
\frac{\partial \varepsilon}{\partial
  t}=-\sum_{l=1}^{L}\frac{w_{l}}{\mu_{R}y_{l}}|\!|\eta_{l}|\!|^{2}\leqq 0.
\label{DerEncont}
\end{equation}
That is, the energy decreases as a function of time, which
expresses the dissipative nature of the problem.

In the discrete case a stability analysis based on energy
techniques permits us to show that the discrete scheme is stable
under the following CFL condition (in homogeneous media and for
both choices $M_h$ and $M_h^1$ for the pressure discretization),
\begin{equation}
\frac{\Delta t^{2}}{4} \frac{\mu_{R}}{\rho} |\!| B_{h} |\!| ^{2}
\left( 1+ \sum_{l=1}^{L} y_{l} \right) \leqq 1 ,
\end{equation}
with $\dsp |\!| B_{h}^{T} B_{h} |\!| \geqq \frac{4}{h^{2}}$ in 1D and $\dsp
|\!| B_{h}^{T} B_{h} |\!| \geqq \frac{8}{h^{2}}$ in 2D.
 Note that these are the usual CFL conditions obtained in the
 non-dissipative case multiplied by $\dsp \left( 1+ \sum_{l=1}^{L} y_{l}
 \right)$.\\
\begin{figure}[htbp]
\label{Dissp}
\begin{minipage}{15cm}
\begin{minipage}{7.5cm}
\centering Dispersion curves for $\pi/4$
\end{minipage}
\begin{minipage}{7.5cm}
\centering Dispersion curves for $0$
\end{minipage}
\end{minipage}~\\
\begin{minipage}{15cm}
\begin{minipage}{7.5cm}
\psfig{figure=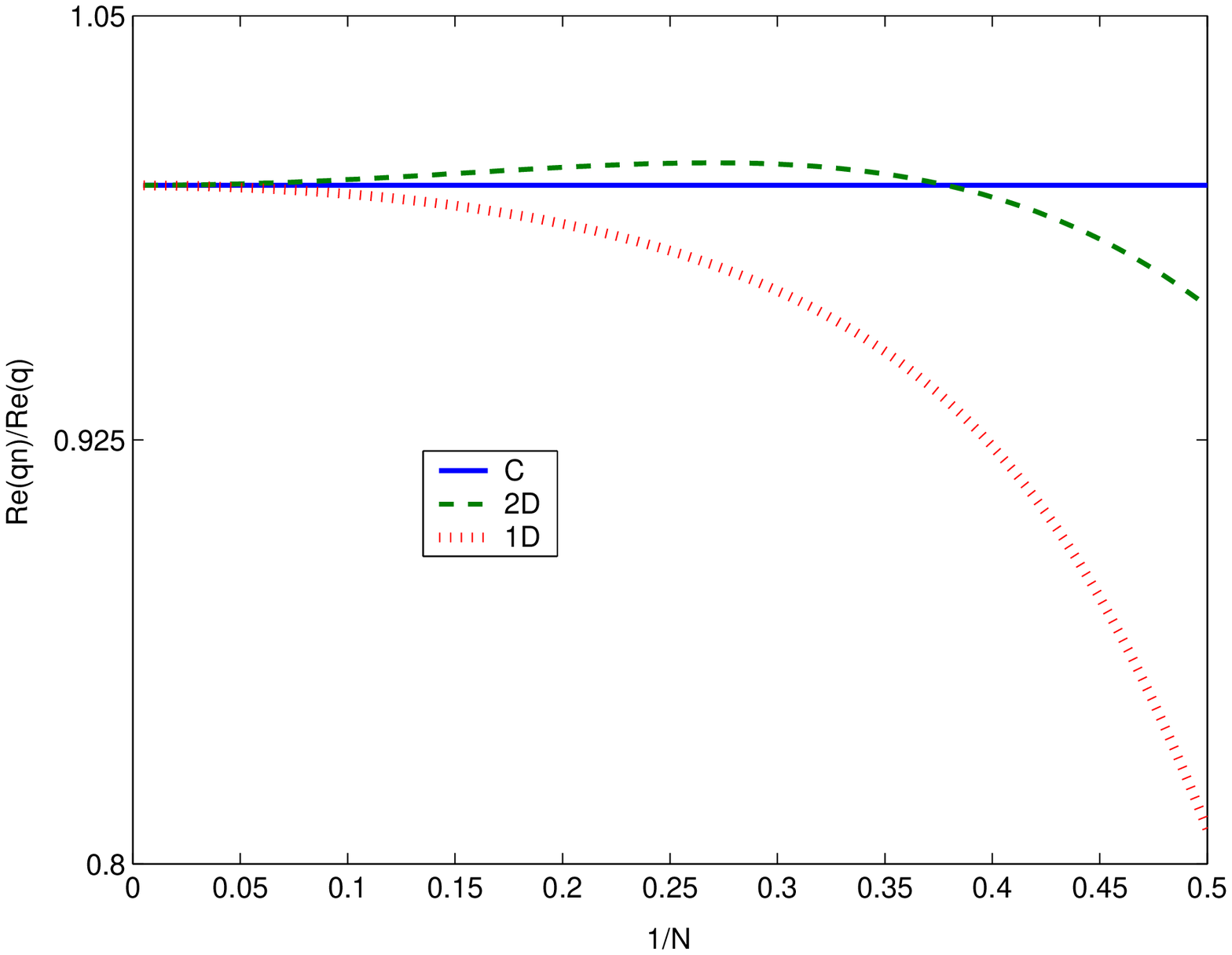,width=7cm}
\end{minipage}
\begin{minipage}{7.5cm}
\psfig{figure=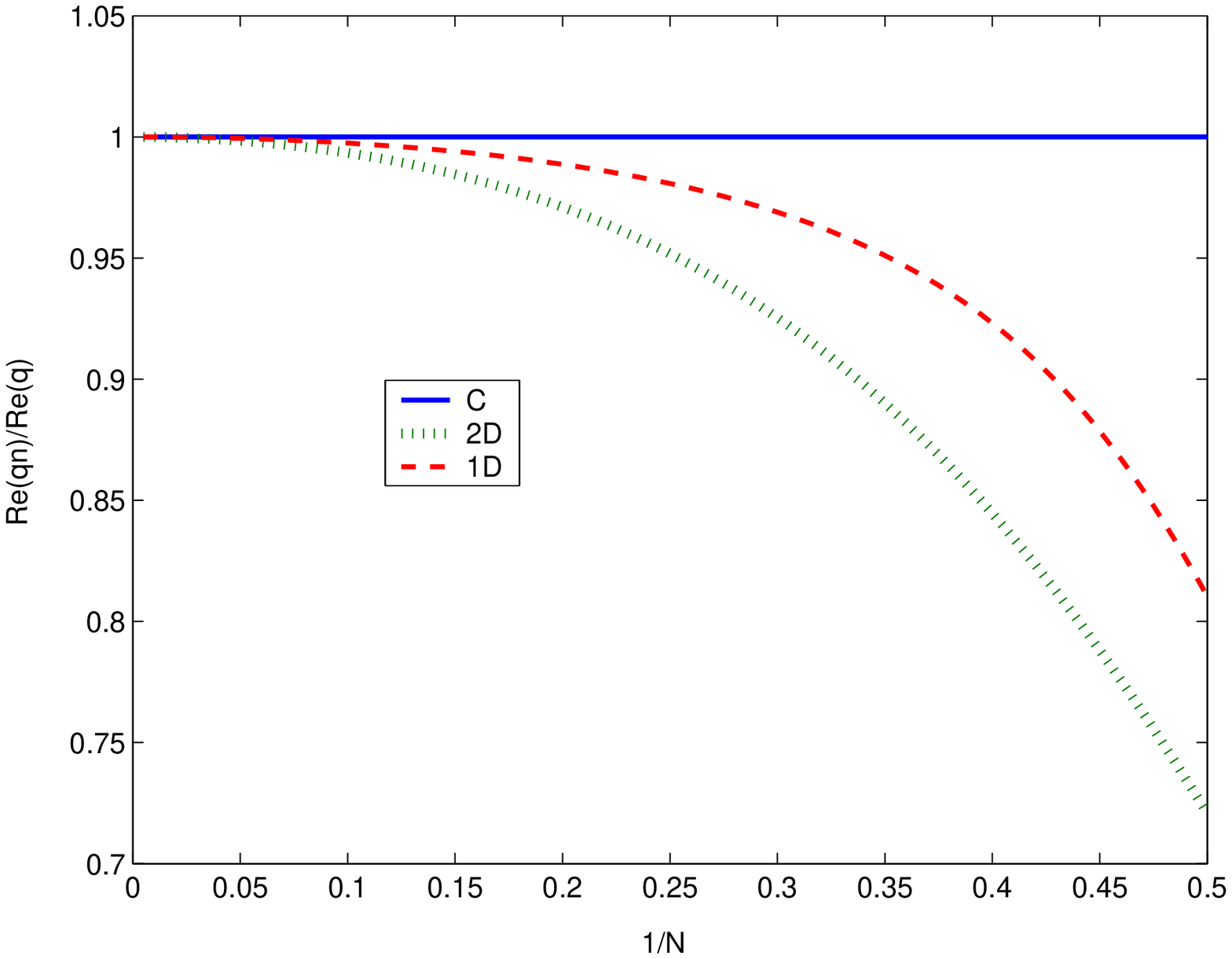,width=7cm}
\end{minipage}
\end{minipage}~\\
\begin{minipage}{15cm}
\begin{minipage}{7.5cm}
\centering Attenuation curves for $\pi/4$
\end{minipage}
\begin{minipage}{7.5cm}
\centering Attenuation curves for $0$
\end{minipage}
\end{minipage}~\\
\begin{minipage}{15cm}
\begin{minipage}{7.5cm}
\psfig{figure=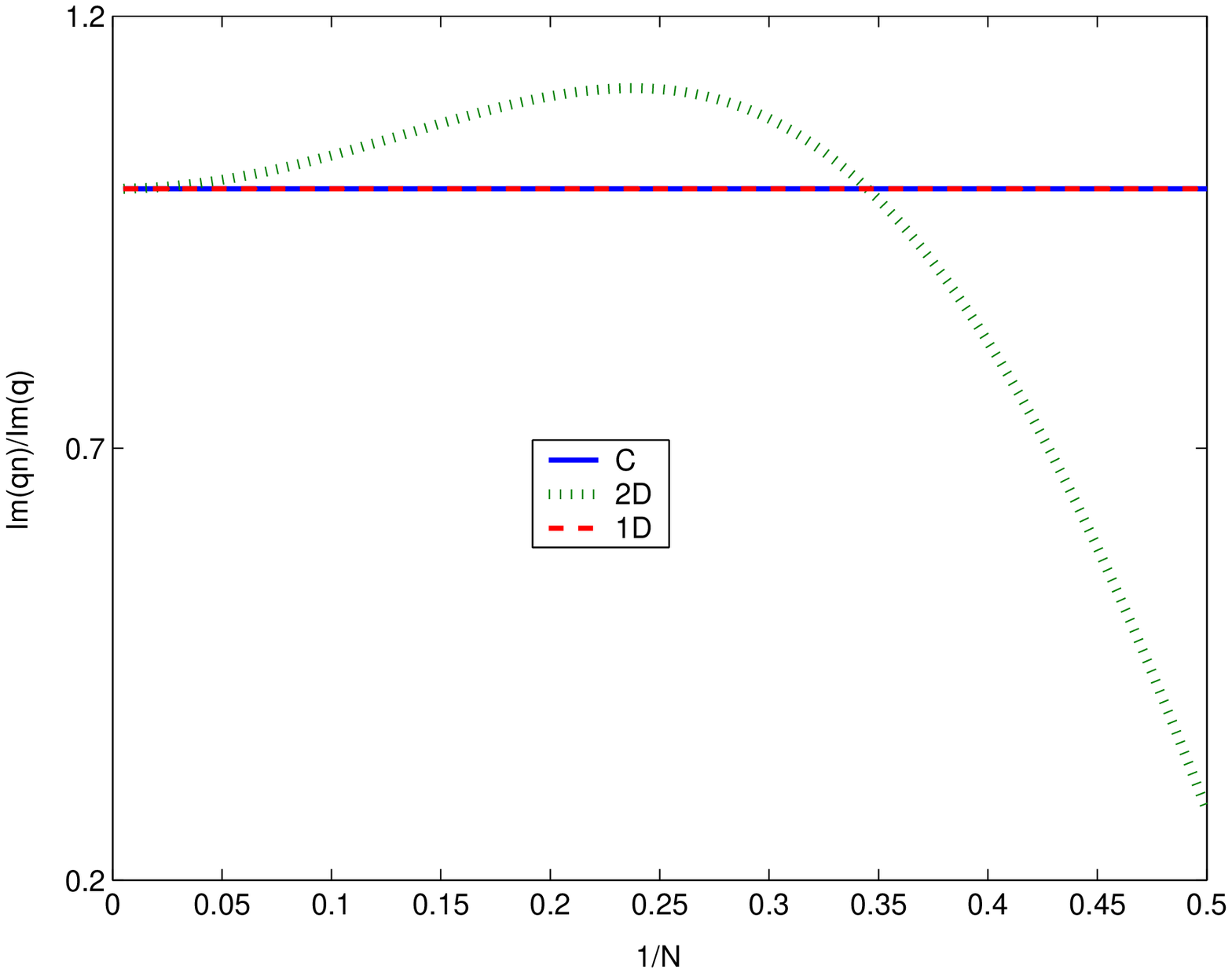,width=7cm}
\end{minipage}
\begin{minipage}{7.5cm}
\psfig{figure=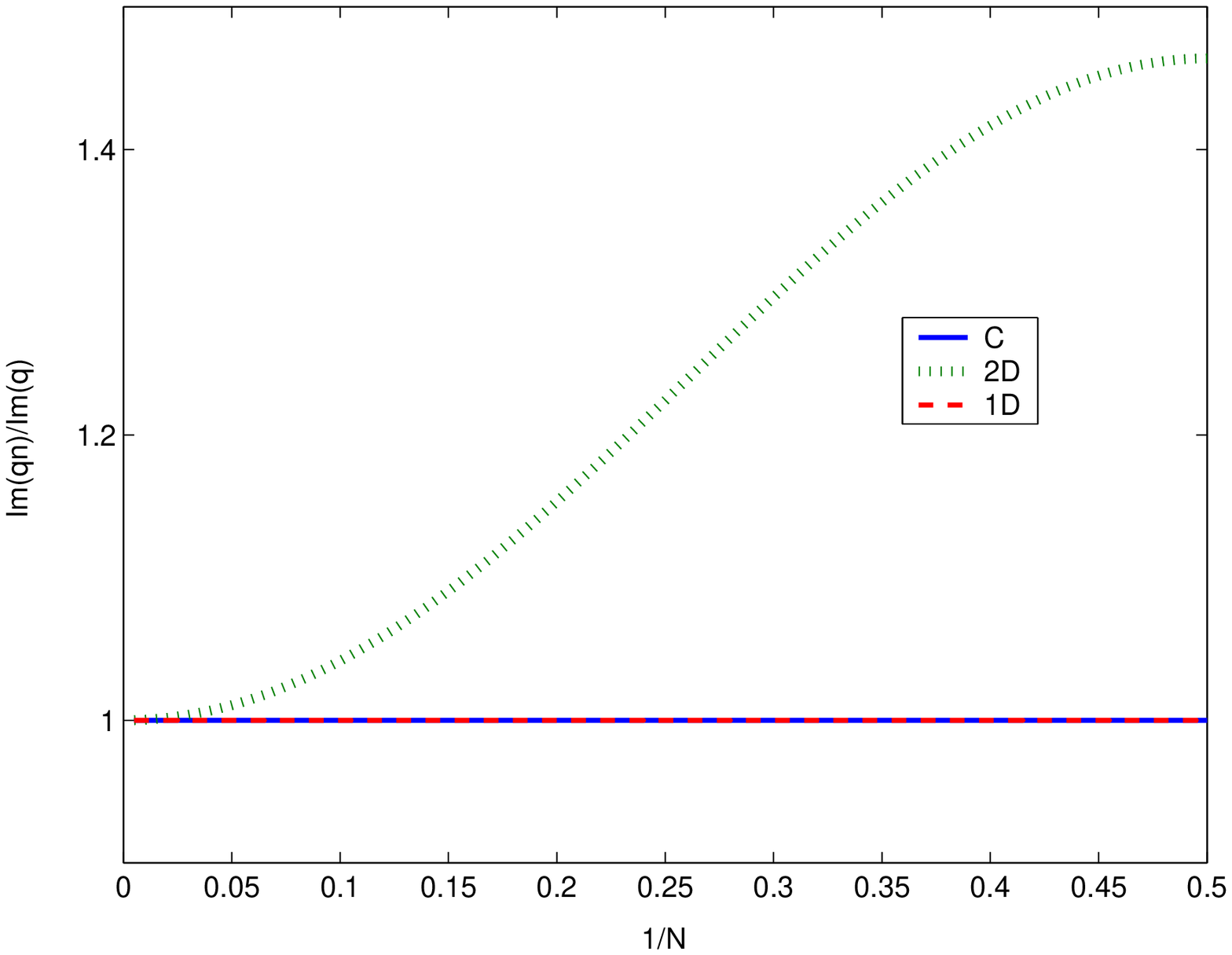,width=7cm}
\end{minipage}
\end{minipage}
\caption{Dispersion and attenuation curves for a plane incident wave,
  with incident angle $0$ or $\pi/4$ for the 2D case.}
\end{figure}
Furthermore, the dispersion relation also presents a similar
aspect. For the the continuous problem we have,
\begin{equation}
\omega^{2}=k^{2} c^{2} \left(1+\sum_{l=1}^{L}\frac{\mathbf{i} \omega
  y_{l}}{\mathbf{i} \omega +\omega_{l}} \right).
\label{Dispcont}
\end{equation}
For the discrete problem in 1D we obtain,
\begin{equation}
\dsp \sin^{2} \left( \frac{\omega \Delta_{t}}{2}
\right)=\frac{\Delta_{t}^2 c^2}{4} \left(\sin^{2} \left( \frac{k
  \Delta_{x}}{2}  \right)  \right) \left(1+\sum_{l=1}^{L}
\frac{2\mathbf{i} y_{l} \tan \left(\frac{\omega \Delta_{t}}{2}
  \right)}{\Delta_{t}\omega_{l}+2 \mathbf{i} \tan \left(\frac{\omega
    \Delta_{t}}{2}  \right) } \right),
\label{Dispdisc1}
\end{equation}
and in 2D we get (for both choices $M_h$ and $M_h^1$ of the
pressure discretization),
\begin{equation}
 \dsp \sin^{2} \left( \frac{\omega \Delta_{t}}{2}
\right)=\frac{\Delta_{t}^2 c^2}{4} \left( \sin^{2} \left( \frac{k_{x}
  \Delta_{x}}{2}  \right) +\sin^{2} \left( \frac{k_{y} \Delta_{y}}{2}
\right) \right) \left(1+\sum_{l=1}^{L} \frac{2\mathbf{i} y_{l} \tan
  \left(\frac{\omega \Delta_{t}}{2}  \right)}{\Delta_{t}\omega_{l}+2
  \mathbf{i} \tan \left(\frac{\omega \Delta_{t}}{2}  \right) } \right).
\label{Dispdisc2}
\end{equation}

In figure \ref{Dissp} we have plotted the
dispersion and attenuation curves as function of $1/N$ ($N$ being
the number of points per wavelength used in the discretization)
for a plane incident wave, whose incident angle is $0$ or
$\pi/4$ for the 2D case. Note, in particular, that the 1D
scheme is no longer exact as it is the case in non-dissipative
media. Depending on the angle of incidence, the 2D scheme may be
more or less dispersive than the 1D one.\\

Demonstration and details of the calculations for the stability
and the dispersion relations for the discrete problem are exposed
in the Appendix \ref{Ap1} and \ref{Ap2}.
\section{The fictitious domain method}{\label{s7}}
To model the free-surface boundary condition on the surface of the
earth we use the fictitious domain method which
has been developed for solving problems involving complex geometries
\cite{B73, GPP94a, GPP94b, GG95, GK98}, and in particular for wave
propagation problems \cite{CJM97, Sylv, Leil, BJT4}.

We follow here the approach proposed in \cite{BJT4}. Consider the
viscoacoustic wave propagation problem in a domain with a complex
geometry such as the one described in Figure \ref{Domfict}. The
initial problem is posed in $\Omega$ with the free-surface
boundary condition, $\mathbf{v} \cdot \mathbf{n} =0$ on $\Gamma$,
\begin{equation}
\left\{ \begin{array}{ll}
\displaystyle \rho \frac{\partial \mathbf{v}}{\partial t}-\nabla p =
\mathbf{f}&\text{in }\Omega, \\[12pt]
\displaystyle \frac{\partial p}{\partial t} - \sum_{l=1}^{n}
\frac{\partial \eta_{l}}{\partial t} = \mu_{R} \Div(\mathbf{v}) &
\text{in } \Omega, \\[12pt]
\displaystyle \frac{\partial \eta_{l}}{\partial t} +
\omega_{l}\eta_{l} = \mu_{R}y_{l} \Div(\mathbf{v}), \forall l & \text{in
} \Omega ,
\\[12pt]
\mathbf{v} \cdot \mathbf{n} = 0, & \mbox{on } \Gamma,\\[6pt]
p=0, &\mbox{on } \Gamma_D.
\label{sys1}
\end{array}
\right.
\end{equation}

\begin{figure}[htbp]
\centering\psfig{figure=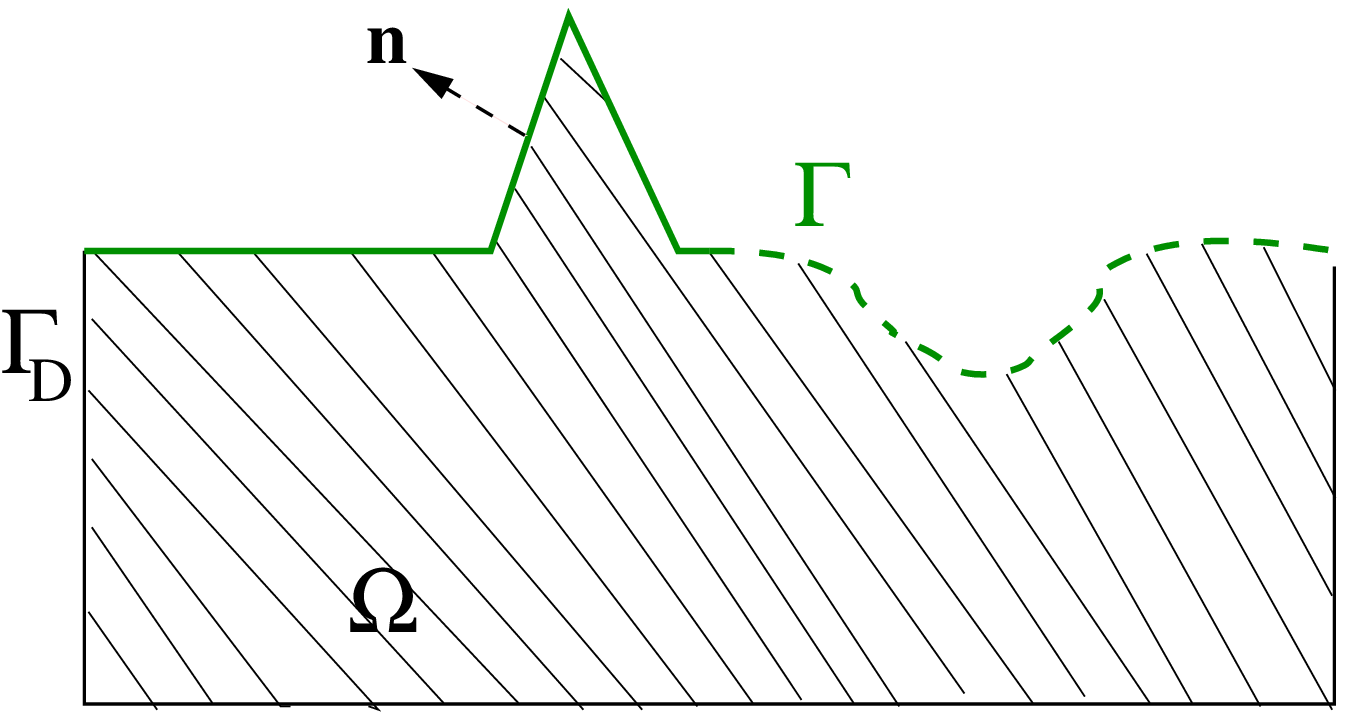,width=7cm}
\psfig{figure=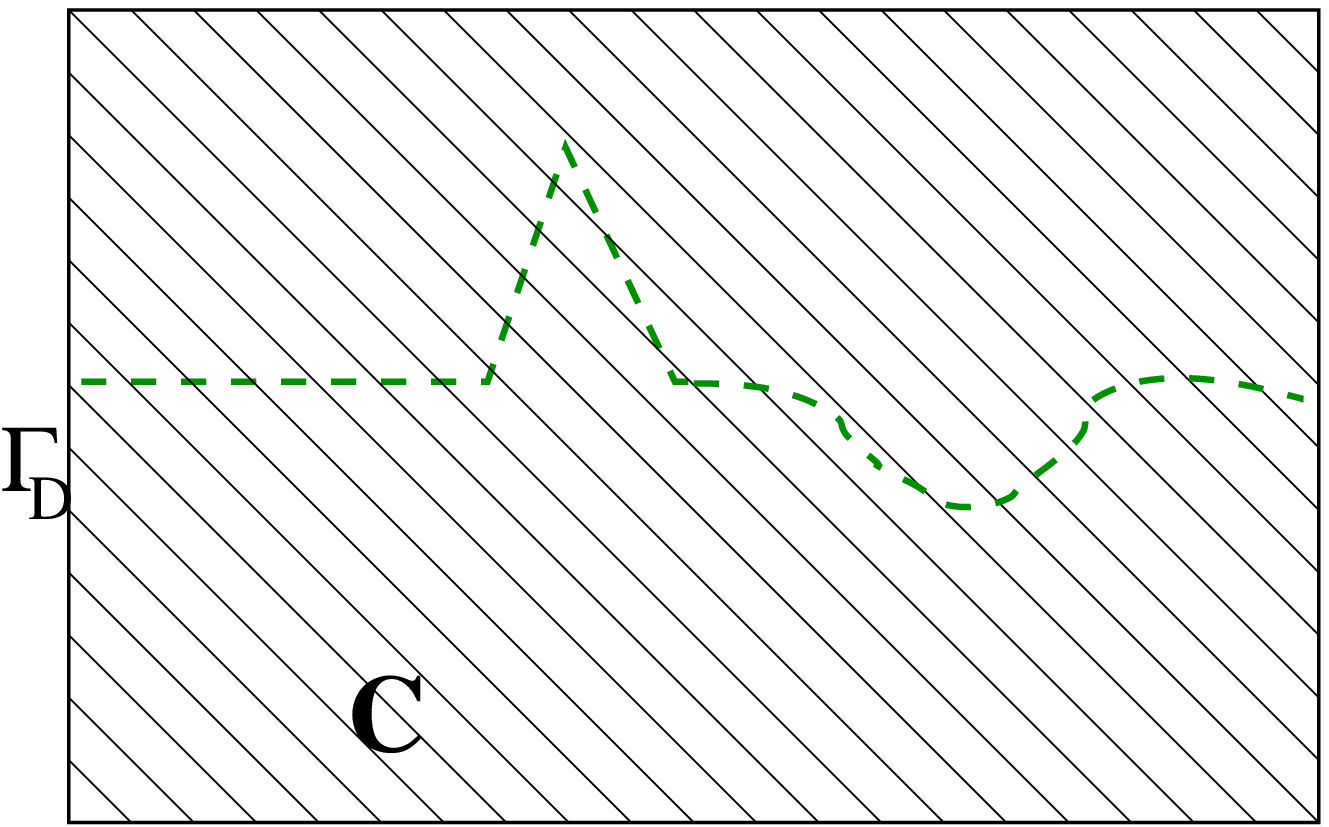,width=7cm}
\caption{Geometry of the problem: on the left the initial domain
  $\Omega$ and on the right the extended domain $C$.}
\label{Domfict}
\end{figure}

The main idea of the fictitious domain method is to extend the
solution to a domain with a simple shape, independent of the
complex geometry of the boundary, and to impose the boundary
condition in a weak way by introducing a Lagrange multiplier.
Following this idea, we extend the solution
$(\mathbf{v},p,\eta_{l})$ by zero in the domain $C$ (which is here
a rectangle, see Figure \ref{Domfict}). We denote
$(\tilde{\mathbf{v}},\tilde{p},\tilde{\eta}_{l})$ the extended
solution and have, $$ \displaystyle \left[ \tilde{\mathbf{v}}
\mathbf{n} \right]_{\Gamma}  =
 0 \Rightarrow  \tilde{\mathbf{v}} \in H(div,C)\text{,   }
 \displaystyle \left[ \tilde{p} \right]_{\Gamma}  \neq 0 \text{,   }\left[
 \tilde{\eta}_{l} \right]_{\Gamma}  \neq 0.
$$

Thus, system (\ref{sys1}) for the extended solution, can be
written (in the distributional sense),

\begin{equation}
\left\{ \begin{array}{ll}
\displaystyle \rho \frac{\partial \tilde{\mathbf{v}}}{\partial
  t}-\nabla \tilde{p} = \mathbf{f} + [\tilde{p}] \mathbf{n}
\delta_{\Gamma}  &\text{in } C, \\[12pt]
\displaystyle \frac{\partial \tilde{p}}{\partial t} - \sum_{l=1}^{n}
\frac{\partial \tilde{\eta}_{l}}{\partial t} = \mu_{R}
\Div(\tilde{\mathbf{v}}) & \text{in } C, \\[12pt]
\displaystyle \frac{\partial \tilde{\eta}_{l}}{\partial t} +
\omega_{l} \tilde{\eta}_{l} = \mu_{R}y_{l} \Div(\mathbf{v}), \forall l
& \text{in } C ,
\\[12pt]
\tilde{\mathbf{v} }\cdot \mathbf{n} = 0, & \mbox{on } \Gamma, \\[6pt]
p=0, &\mbox{on } \Gamma_D.
\label{sys2}
\end{array}
\right.
\end{equation}

In (\ref{sys2}) we have two types of unknowns, the extended
unknowns, defined in the simple shape domain $C$ and the auxiliary
variable $[\tilde{p}]$, defined on the boundary $\Gamma$. We
introduce $\lambda=[\tilde{p}]$ as a new unknown defined on
$\Gamma$. This unknown can be interpreted as a Lagrange multiplier
associated with the boundary condition on $\Gamma$. The
variational formulation of the problem can then be written as
follows,

$$\left\{
\begin{array}{ll}
\text{Find } (\mathbf{v},p,H,\lambda) :  ]0,T[  \longmapsto  H(\Div;C)
  \times L^2(C) \times (L^2(C))^{L} \times H^{1/2}(\Gamma) \text{
  s.t. } & \\[6pt]
\displaystyle \frac{d}{dt} (\rho \mathbf{v},\mathbf{w}) +
b(\mathbf{w},p)  - b_{\Gamma}(\lambda,\mathbf{w}) =
    (\mathbf{f},\mathbf{w}) , & \displaystyle
    \forall {\mathbf w} \in H(\Div;C) , \\ [12pt]
\displaystyle \frac{d}{dt} (\frac{1}{\mu_{R}} p,q) - \sum_{l=1}^{L}
    \frac{d}{dt}(\frac{1}{\mu_{R}} \eta_{l},q) - b({\mathbf{v}},q)=0,
    & \displaystyle \forall q \in L^2(C) , \\ [12pt]
\displaystyle \frac{d}{dt} (\frac{1}{\mu_{R} y_{l}} \eta_l,q)
+ (\frac{\omega_{l}}{\mu_{R} y_{l}} \eta_{l},q) -
    b({\mathbf{v}},q)=0,\ \forall l,
& \displaystyle \forall q \in L^2(C) , \\[12pt]
\dsp b_{\Gamma}(\mu,\mathbf{v}) =0, & \forall \mu \in H^{1/2}(\Gamma),
\end{array}
\right.
$$
where

$$\begin{array}{ll}
\displaystyle b_{\Gamma}(\mu,\mathbf{w}) = \int_{\Gamma}  \mu \
\mathbf{w}\cdot \mathbf{n} ds,
& \displaystyle  \forall (\mu,\mathbf{w}) \in
H^{1/2}(\Gamma) \times  H(\Div;C).
\end{array}
$$

For the discretization of this problem we consider a structured
volume mesh ${\cal T}_{h}$ on $C$, and an irregular surface mesh
${\cal G}_{h_s}$ on $\Gamma$. The main advantage of this
formulation is that the mesh for computing the extended functions
can now be regular while the surface mesh is irregular and permits
a good and efficient approximation of the geometry (see Figure
\ref{Domfict2}).

\begin{figure}[htbp]
\centering\psfig{figure=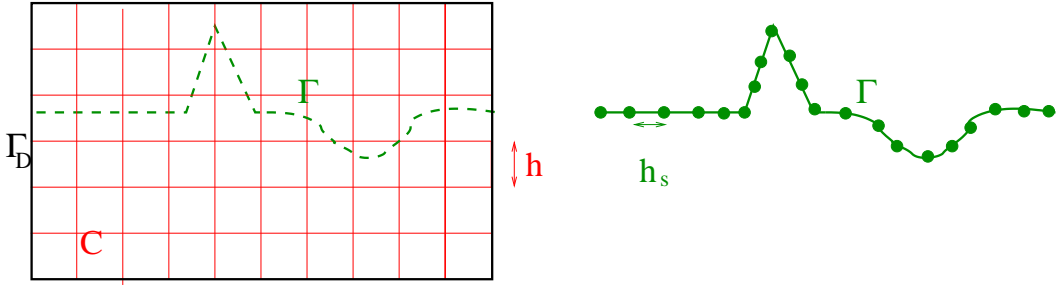,width=14cm}
\caption{The two meshes used in the fictitious domain method: a
  structured volume mesh on the domain $C$ and an irregular surface
  mesh on $\Gamma$.}
\label{Domfict2}
\end{figure}

To discretize the volume unknowns $(\mathbf{v},p,H)$ we use the
finite element method described in section \ref{s5} while for the
Lagrange multiplier $\lambda$ we use piecewise linear continuous
functions on ${\cal G}_{h_s}$, i.e., the approximation space is,
$$  G_{h_s}=\left\{ {\mu}_{h_s} \in H^{1/2}(\Gamma) \left/ \right.
\forall  S (\mbox{segment}) \in  {\cal G}_{h_s},  {\mu}_{h_s}
\left|_{S}
    \right. \in {P_1}(S) \right\}.$$

To simplify the presentation, we considered in system (\ref{sys1}) the
homogeneous Dirichlet boundary condition on the boundary
$\Gamma_D$. When the domain is infinite we use the perfectly matched
layer  model which will be described in the following section.

\section{The PML method}{\label{s6}
The perfectly Matched Layer model was introduced by B\'erenger
\cite{ber94, ber96} for  Maxwell's equations and is  now the most
widely-used method for simulating wave propagation in unbounded
domains. The reader can refer to \cite{zhao,
  teix2, Petr} for electromagnetic waves, to \cite{BFP} for anisotropic
acoustic waves and to \cite{CT, BFP} for
elastic waves. The popularity of this model is due to its simplicity
and efficiency. Its most astonishing property is that for the
continuous problem the reflection coefficient at the interface between
the layer and the free medium is zero for all frequencies and angles
of incidence.
\begin{figure}[htbp]
\centering\psfig{figure=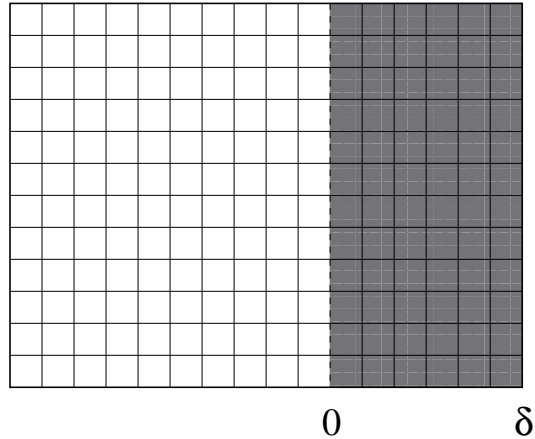,width=7cm}
\caption{PML in the $x$-direction: the physical medium is on the left
  and the absorbing medium is a layer of width $\delta$.}
\label{Pml}
\end{figure}
To derive the PML for the viscoacoustic system (\ref{systime}) we
follow the approach proposed in \cite{CT} which applies to any
first-order linear hyperbolic system. We present here the construction of
a PML in the $x$-direction (see Figure \ref{Pml}). Deriving then the PML
for the other boundaries and the corners of the
computational domain is a straightforward application of the same
technique.

Following \cite{CT} we construct the PML model in two steps: (i)
We split the solution $(\mathbf{v},p,\{\eta_l\} )$ into two parts
$ (\mathbf{v}^{\parallel},p^{\parallel},\{\eta_l^{\parallel}\})$
and $(\mathbf{v}^{\perp},p^{\perp},\{\eta_l^{\perp}\})$, with the
parallel part being associated with the derivatives in the
$y$-direction (direction parallel to the interface between the PML
and the physical medium), and the orthogonal part associated with
those in the $x$-direction. (ii) We introduce  damping only on the
orthogonal component of the solution.

When applying the splitting step to (\ref{systime}) by remarking
that $\mathbf{v}^{\parallel} = (0,v_y)$ and  $\mathbf{v}^{\perp} =
(v_x,0)$, we obtain,
\begin{equation}
\left\{ \begin{array}{l}
\displaystyle \rho \frac{\partial v_y}{\partial
  t}=\frac{\partial p}{\partial y}\\[12pt]
\displaystyle \frac{\partial p^{\parallel} }{\partial
  t}-\sum_{l=1}^{L}\frac{\partial \eta_{l}^{\parallel}}{\partial
  t}=\mu_{R}\frac{\partial v_y}{\partial y} \\[12pt]
\displaystyle \frac{\partial \eta_{l}^{\parallel}}{\partial
  t}+\omega_{l}\eta_{l}^{\parallel}=\mu_{R} y_{l} \frac{\partial
  v_y}{\partial y},
\label{systpml1}
\end{array}
\right.
\end{equation}
and
\begin{equation}
\left\{ \begin{array}{l}
\displaystyle \rho \frac{\partial v_x}{\partial
  t}=\frac{\partial p}{\partial x}\\[12pt]
\displaystyle \frac{\partial p^{\perp} }{\partial
  t}-\sum_{l=1}^{L}\frac{\partial \eta_{l}^{\perp}}{\partial
  t}=\mu_{R}\frac{\partial v_x}{\partial x} \\[12pt]
\displaystyle \frac{\partial \eta_{l}^{\perp}}{\partial
  t}+\omega_{l}\eta_{l}^{\perp}=\mu_{R} y_{l} \frac{\partial
  v_x}{\partial x},
\label{systpml2}
\end{array}
\right.
\end{equation}
with
\begin{equation}
\left\{ \begin{array}{l}
\displaystyle p=p^{\parallel}+p^{\perp}\\
\displaystyle \eta_{l}=\eta_{l}^{\parallel}+\eta_{l}^{\perp}\text{,
}\forall l.\\
\end{array}
\right.
\end{equation}
To apply the damping on the orthogonal components it is simpler to
consider system (\ref{systpml2}) in the frequency domain. Then the
PML consists in replacing the $x$-derivatives $\partial_x$ by $\dsp \frac{i
  \omega}{i \omega+d(x)} \partial_x$ (cf. \cite{CT}). Following this
approach, system (\ref{systpml2}) in the frequency domain becomes,
\begin{equation}
\left\{ \begin{array}{ll}
(i)&\displaystyle \rho (i \omega + d(x)) \partial v_x = \frac{\partial
    p}{\partial x}\\[12pt]
(ii)&\displaystyle (i \omega + d(x)) p^{\perp}
  -\sum_{l=1}^{L} (i \omega + d(x))  \eta_{l}^{\perp} =
\mu_{R}\frac{\partial v_x}{\partial x} \\  [12pt]
(iii)&\displaystyle (i \omega) (i \omega + d(x))  \eta_{l}^{\perp}
 + (i \omega + d(x))  \omega_{l} \eta_{l}^{\perp}= (i \omega ) \mu_{R} y_{l} \frac{\partial
  v_x}{\partial x},
\label{systpml2b}
\end{array}
\right.
\end{equation}
where $d(x)$ is the damping parameter which is equal to zero in the
physical medium and non-negative in the absorbing medium.

We now introduce new variables $\widetilde{\eta}_l$ defined by,
\begin{equation}
\label{eta_t}
i \omega  \widetilde{\eta}_l = (i \omega + d(x)) \eta_{l}^{\perp},\
\forall l,
\end{equation}
or equivalently in time domain,
$$\dsp \frac{\partial \widetilde{\eta}_l}{\partial t} = \frac{\partial
  {\eta}^{\perp}_l}{\partial t} + d(x)  {\eta}^{\perp}_l,\ \forall l.$$
Using (\ref{eta_t}) in (\ref{systpml2b}) and going in the time domain
we get,
\begin{equation}
\left\{ \begin{array}{l}
\displaystyle \rho \frac{\partial v_x }{\partial t}+ \rho
d(x)  v_x = \frac{\partial
    p}{\partial x}\\ [12pt]
\displaystyle \frac{\partial p^{\perp} }{\partial t} + d(x) p^{\perp}
  -\sum_{l=1}^{L} \frac{\partial \widetilde{\eta}_{l}}{\partial t} =
\mu_{R}\frac{\partial v_x}{\partial x} \\   [12pt]
\displaystyle \frac{\partial \widetilde{\eta}_{l}}{\partial t}
 +  \omega_{l} \widetilde{\eta}_{l} =  \mu_{R} y_{l} \frac{\partial
  v_x}{\partial x}.
\label{systpml2c}
\end{array}
\right.
\end{equation}
The final system of equations for the PML is (\ref{systpml2c})
together with (\ref{systpml1}), with $p$ being defined by
$p=p^{\parallel}+p^{\perp}$. Note that the memory variables
$\eta_l$ do not appear, and only the component
$\eta_l^{\parallel}$ and the variables $\widetilde{\eta}_l$ do
appear, in this system.

Using a plane wave analysis, it can be shown (cf. \cite{CT}) that
this model generates no reflection at the interface between the
physical and the absorbing medium and that the wave decreases
exponentially inside the layer.  This property allows the use of a
very high damping parameter inside the layer, and consequently of
a small layer width, while achieving a near-perfect absorption of
the waves. Note that for a finite-length absorbing layer there is
some reflection due to the outer boundary of the PML.

\begin{remark} To discretize the PML we use the same scheme as for the
interior domain.\end{remark}

\begin{remark}
The damping $d(x)$ is zero in the physical domain and non negative in
the absorbing medium. In the numerical simulations it is defined as in
\cite{CT},
\begin{equation}
d(x)=\left\{\barr{ll}
0 & \text{for } x < 0\\
log \left(\frac{1}{R} \right) \frac{(n+1)
  \sqrt{\frac{\mu_R}{\rho}}}{2 \delta} (\frac{x}{\delta})^n
&\text{for } x \ge 0
\earr\right.
\end{equation}
where $R$ is the theoretical reflection coefficient, $\delta$ the
width of the PML and $n=4$.\\

In practice, we take $R=5.0 \ 10^{-7}$, and $\delta \approx
30\Delta_{x}$ (depending on the wavelength).
\end{remark}

\section{Numerical results}{\label{s8}
\subsection{Scattering from a circular cylinder}{\label{s81}
In order to validate the proposed numerical method we consider in
this section the canonical problem of a plane wave (Ricker
wavelet) striking a viscoacoustic homogeneous circular cylinder.
The geometry of the problem is displayed in Figure \ref{Diffcyl}.
A homogeneous viscoacoustic circular cylinder of radius $a$
(domain $\Omega_2$) is surrounded by a homogeneous non-dissipative
medium (domain $\Omega_1$). We denote by $\Gamma_1$ the interface
between the two domains $\Omega_1$ and $\Omega_2$. The physical
characteristics of the media are $\rho_1=1000$Kgr/m$^3$,
$c_1=3050$m/s, $Q_1 = +\infty$ in $\Omega_1$ and
$\rho_2=1800$Kgr/m$^3$, $c_2=3050$m/s and $Q_2= 30$ in $\Omega_2$.
The source function used in this example is given by
(\ref{ricker}) with $f_0=2.5$Hz. For this problem, the solution
can be computed by an analytical method  described in what
follows.

\begin{figure}[htbp]
\centering\psfig{figure=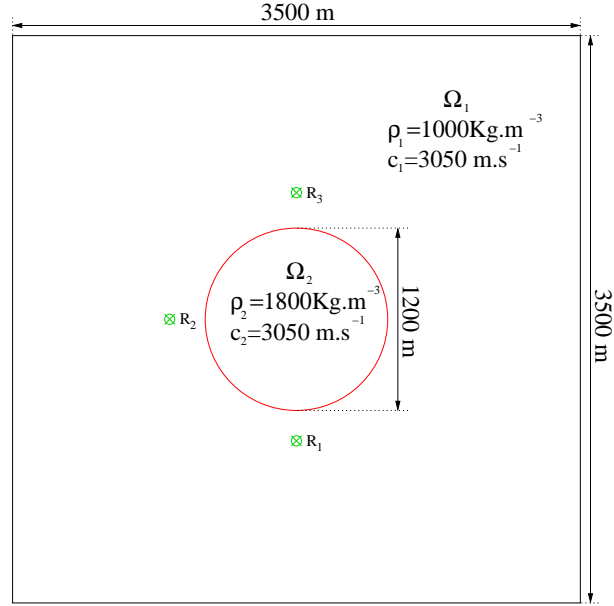,width=8cm}
\caption{The geometry of the problem: a homogeneous viscoacoustic
  cylinder of radius $a$ (domain $\Omega_2$) embedded in a
  non-dissipative homogeneous medium (domain $\Omega_1$).}
\label{Diffcyl}
\end{figure}
Consider the following incident plane wave (with incident angle $0 \le
\theta^{i} < 2 \pi$),
\begin{equation}
p_{1}^{i}(\mathbf{x})=
A_{0}^{i} \sum_{n=-\infty}^{n=\infty} e^{-\mathbf{i}n
  (\theta^{i}+\frac{\pi}{2})} J_{n}(k_{0}r) {\exp(\mathbf{i}n \theta)} \text{;
} \forall \mathbf{x}=(r\cos(\theta),r\sin(\theta)) \in \Omega_{1}.
\end{equation}
Using the partial wave expansion we can express the solutions $p_j
\in \Omega_j, j=1, 2$ in the following form,
\begin{equation}
\begin{array}{ll}
\displaystyle
p_{1}(\mathbf{x})=p_{1}^{i}(\mathbf{x})+\sum_{n=-\infty}^{n=\infty}
a_{1n} H_{n}^{(1)}(k_{1}r) {\exp (\mathbf{i} n \theta)},
&
\forall \mathbf{x} \in \Omega_{1}\\[6pt]
\displaystyle p_{2}(\mathbf{x})
=\sum_{n=-\infty}^{n=\infty}a_{1n}J_{n}(k_{2}(\omega)r){\exp(\mathbf{i}n\theta)},
&
\forall \mathbf{x} \in \Omega_{2},
\end{array}
\label{Diffcyl1}
\end{equation}
with $H_n^{(1)}$ the first-kind Hankel function of
order $n$ , $J_n$ the Bessel function of order $n$ and
where the wave number in $\Omega_{2}$ is computed by (\ref{Kjar}),
\begin{equation}
k_{2}(\mathbf{x},\omega)=k_{2}(\omega)=\frac{\omega}{c_{ref}} \left(
\frac{\mathbf{i}\omega}{\omega_{ref}}\right)^{-\frac{1}{\pi}atan(\frac{1}{Q_{2}})}.
\end{equation}
To compute the coefficients $a_{1n}$ and $b_{1n}$ we introduce the
expressions for $p_1$ and $p_2$, i.e., equation (\ref{Diffcyl1})
in the transmission boundary conditions on $\Gamma_1$ (continuity
of the pressure and the normal component of velocity). After
projecting the resulting system onto the Fourier basis $\dsp
\left( \frac{1}{2\pi} {\exp(-\mathbf{i}m \theta)} \text{;  } m\in
\mathbf{Z} \right)$ we obtain,
\begin{equation}
\begin{array}{l}
\displaystyle a_{1n}= \frac{
\gamma_{1}\dot{J}_{n}(\chi_{1}) J_{n}(\chi_{2}) -
\gamma_{2} J_{n}(\chi_{1})\dot{J}_{n}(\chi_{2})}
{\gamma_{1} \dot{H}_{n}^{(1)}(\chi_{1}) J_{n}(\chi_{2})
+ \gamma_{2}\dot{J}_{n}(\chi_{2}) H_{n}^{(1)}(\chi_{1})} A_{0}^{i}e^{-\mathbf{i} n
  (\theta^{i}+\frac{\pi}{2})}\\[12pt]
\displaystyle b_{2n}= \frac{\gamma_{1} \dot{J}_{n}(\chi_{1})
  H_{n}^{(1)}(\chi_{1})-\gamma_{1} J_{n}(\chi_{1})
  \dot{H}_{n}^{(1)}(\chi_{1})}{\gamma_{2} \dot{J}_{n}(\chi_{2})
  H_{n}^{(1)}(\chi_{1})+\gamma_{1} \dot{H}_{n}^{(1)}(\chi_{1})
  J_{n}(\chi_{2})}A_{0}^{i}e^{-\mathbf{i} n
  (\theta^{i}+\frac{\pi}{2})},
\end{array}
\label{Diffcoef}
\end{equation}
with $\dot{Z}_{n}(z)=\frac{dZ_{n}(z)}{dz}$, $\chi_{j}=k_{j}a$, and
$\gamma_{j}=\frac{k_{j}}{\rho_{j}}$. The insertion of these
expressions into (\ref{Diffcyl1}) gives the final solution of the
problem \cite{Wirgin}. Comparison of results between the
analytical and the numerical solution are displayed in Figure
\ref{Diffcylr1} where we can see that good agreement is obtained
between the two.

\begin{figure}[htbp]
\centering \begin{minipage}{16cm}
\begin{minipage}{5cm}
\hspace*{0.8cm} Solution at point $R_1$
\end{minipage}
\begin{minipage}{5cm}
\hspace*{0.8cm} Solution at point $R_2$
\end{minipage}
\begin{minipage}{5cm}
\hspace*{0.8cm} Solution at point $R_3$
\end{minipage}
\end{minipage}~\\
\begin{minipage}{16cm}
\psfig{figure=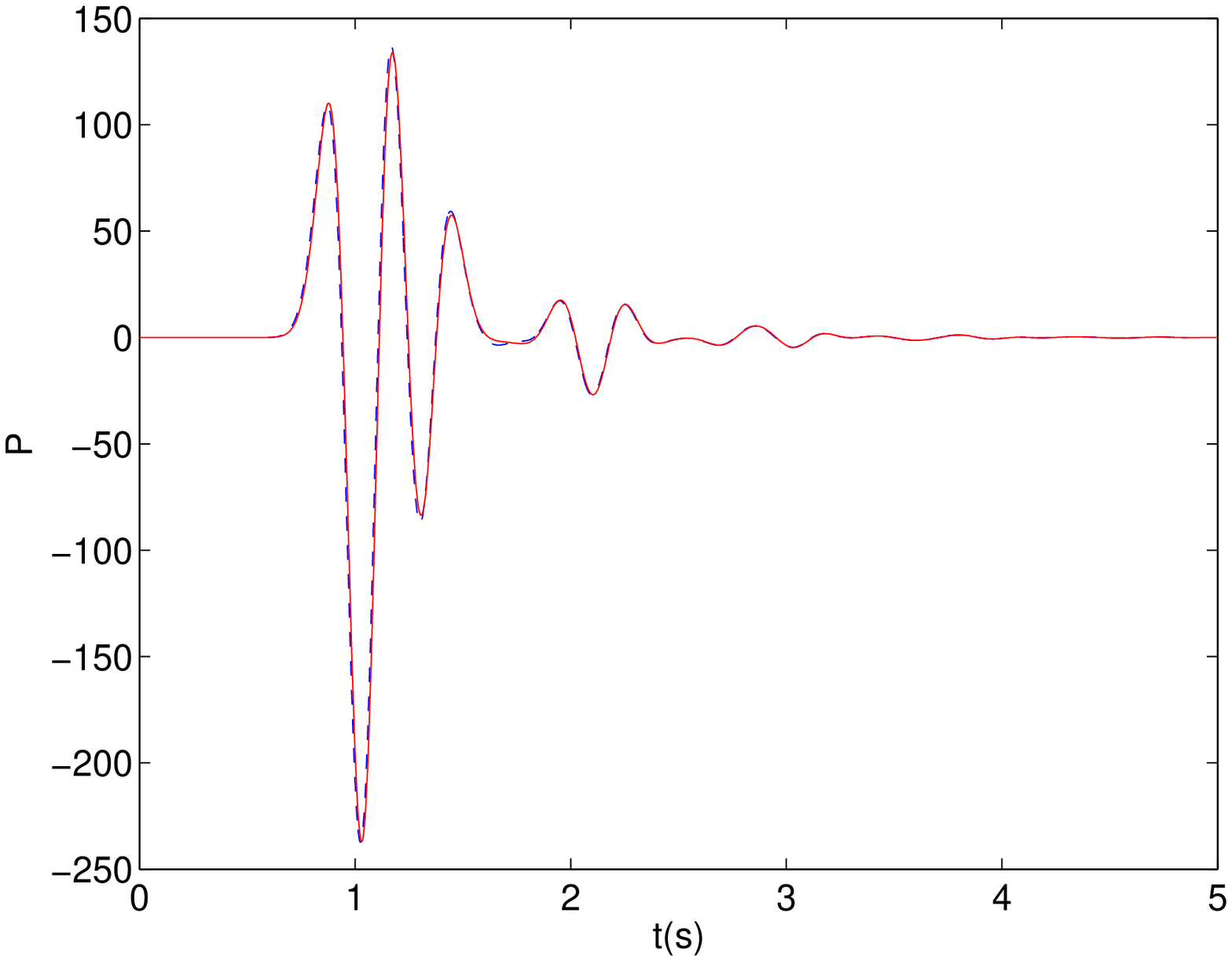,width=5cm}
\psfig{figure=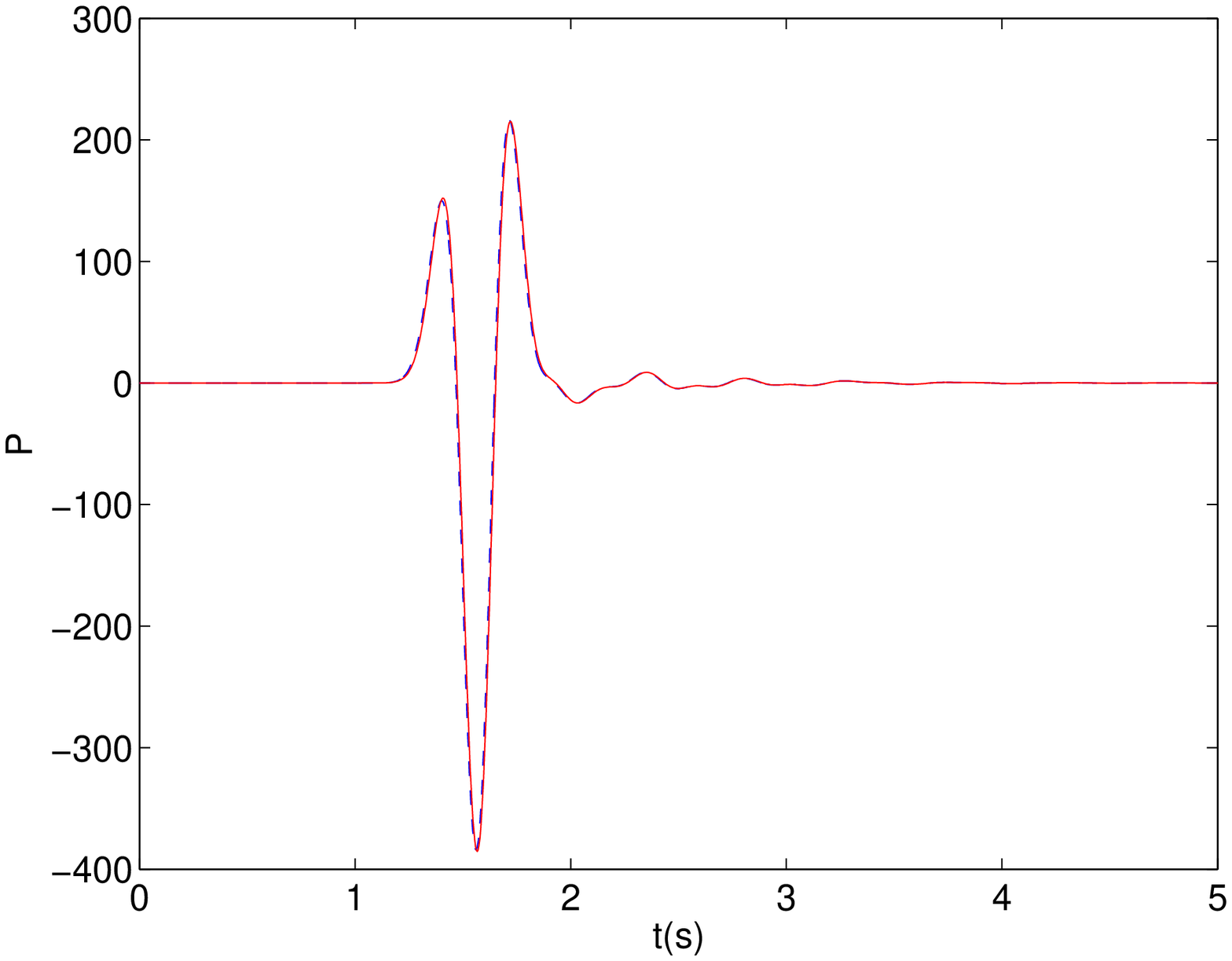,width=5cm}
\psfig{figure=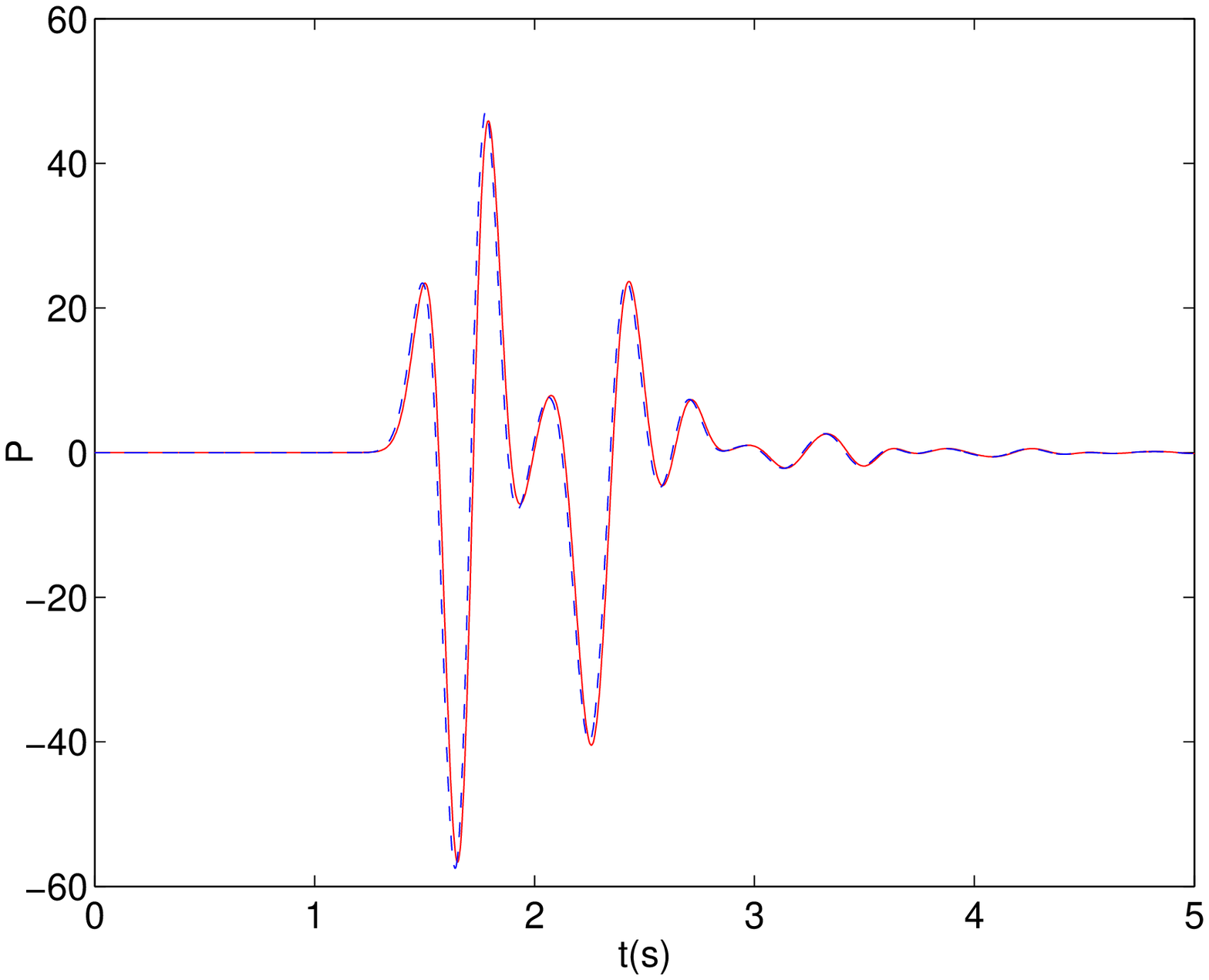,width=5cm}
\end{minipage}
\caption{Comparison between the analytical solution (dashed line) and the
  numerical solution (continuous line) at different observation points $R_1$,
  $R_2$ and $R_3$. The location of the observation points is
  illustrated in Figure \ref{Diffcyl}. In the figures the x-axis is
  time (in s) and the y-axis is the pressure field.}
\label{Diffcylr1}
\end{figure}
In the numerical simulation, we assume that the problem is posed in the whole
space and to solve it, we couple system (\ref{TotDi})
with the perfectly matched absorbing layer model (PML).
\subsection{Simulation of the response to an incident cylindrical wave
  of a dike on a flexible foundation embedded in a half-space}
\label{s82} To illustrate the efficiency of the method we model in
this section the response to an incident cylindrical wave of a
dike on a flexible foundation embedded in a half-space. This
particular problem was considered in \cite{Todo2001} where it was
solved using an expansion of the solution in cylindrical wave
functions in the case of non-dissipative media. In \cite{Todo2001}
the authors studied this problem for different material parameters
in order to determine how stiff the foundation should be relative
to the soil for the rigid foundation assumption in soil-structure
interaction models to be valid. They concluded that a foundation
with the same mass density as the soil but 50 times larger shear
modulus behaves in rigid manner for this problem. However, for
ratios of shear moduli less than 16, the rigid foundation
assumption is not valid. In this case, soil-structure interaction
models with a rigid-foundation assumption will not model the
differential motion of the ground and may underestimate the
stresses in the structure (cf. \cite{Todo2001}).  We consider here
a ratio of shear moduli equal to 4. Soil-structure interaction is
taken into account owing to the fact that we discretize the
continuous problem.

\begin{figure}[htbp]
\centering\psfig{figure=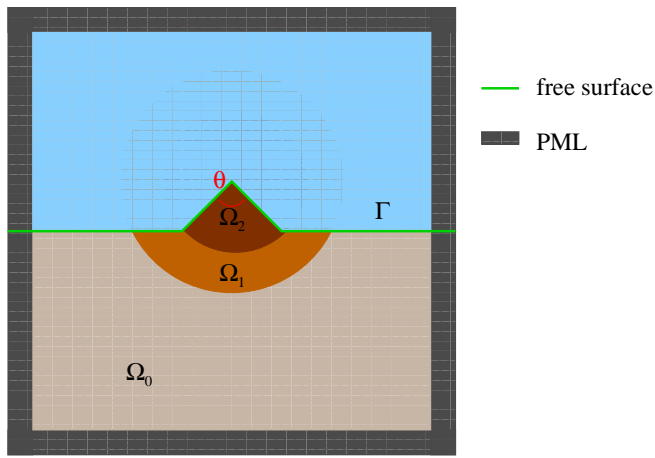,width=9cm}
\caption{The geometry of the problem: a dike
  on a flexible foundation embedded in a half-space.}
\label{Mont1}
\end{figure}

The geometry of the problem is illustrated in Figure \ref{Mont1},
where $\Gamma$ denotes the free surface, $\Omega_0$  the hard
bedrock, $\Omega_1$ the flexible foundation and $\Omega_2$ the
dike. The physical parameters used in the simulation are
$\rho_{0}=1000$Kg/m~$^{3}$, $c_{0}=1450$m/s, $Q_{0}=+\infty$ in
the bedrock, $\rho_{1}=1000$Kg/m~$^{3}$, $c_{1}=2900$m/s,
$Q_{1}=30$ in $\Omega_1$ and $\rho_{2}=250$Kg/m~$^{3}$,
$c_{2}=725$m/s, $Q_{2}=100$ in the dike. The angle $\theta$ is
equal to ${\pi}/{2}$.

\begin{figure}[htbp]
\begin{minipage}{15cm}
\begin{minipage}{7cm}
\hspace*{3.2cm} $t=2$s
\end{minipage}
\begin{minipage}{7cm}
\hspace*{3.2cm} $t=3$s
\end{minipage}
\end{minipage}~\\
\begin{minipage}{15cm}
\begin{minipage}{7cm}
\centering\psfig{figure=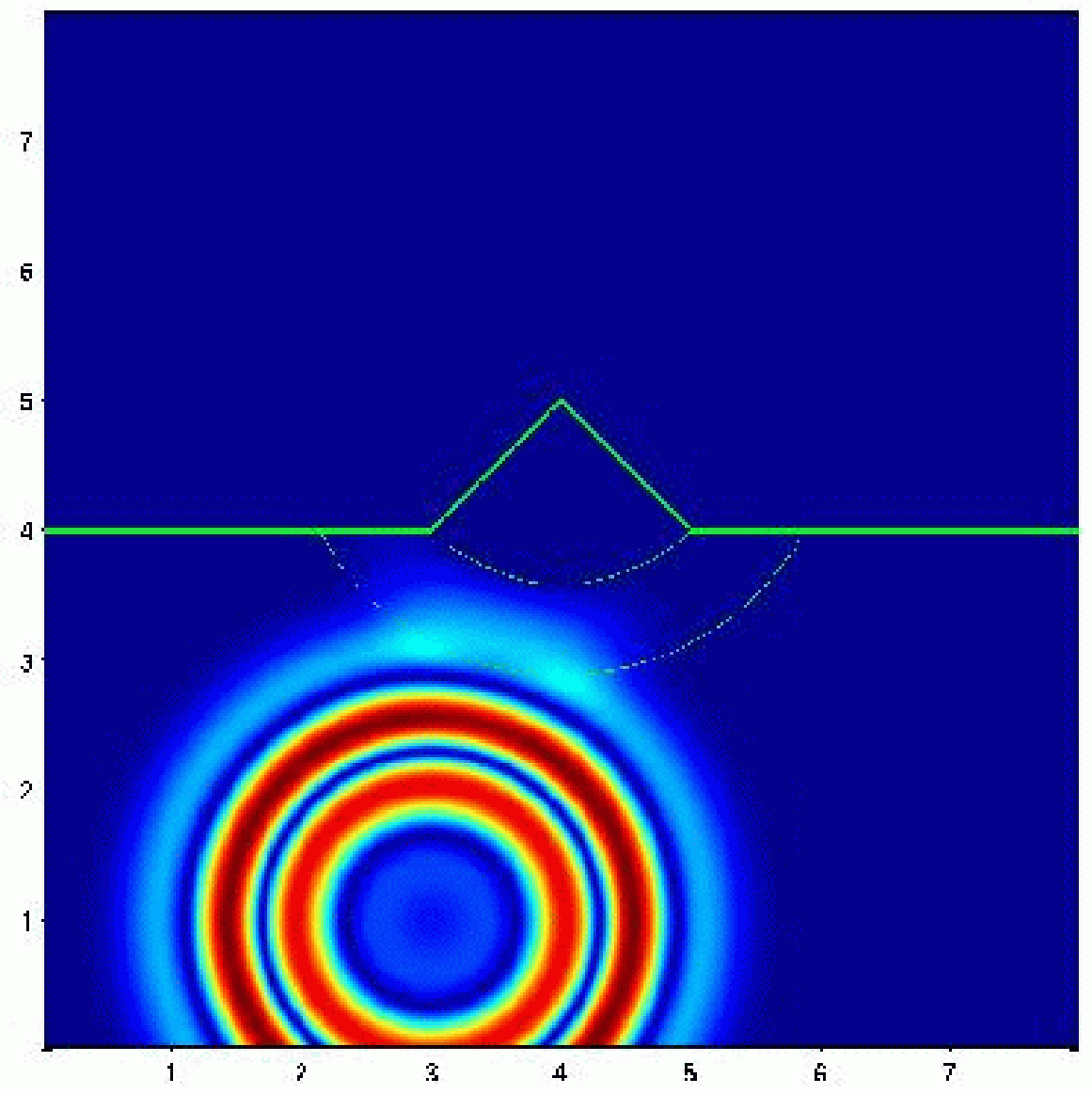,width=7cm}
\end{minipage}
\begin{minipage}{7cm}
\centering\psfig{figure=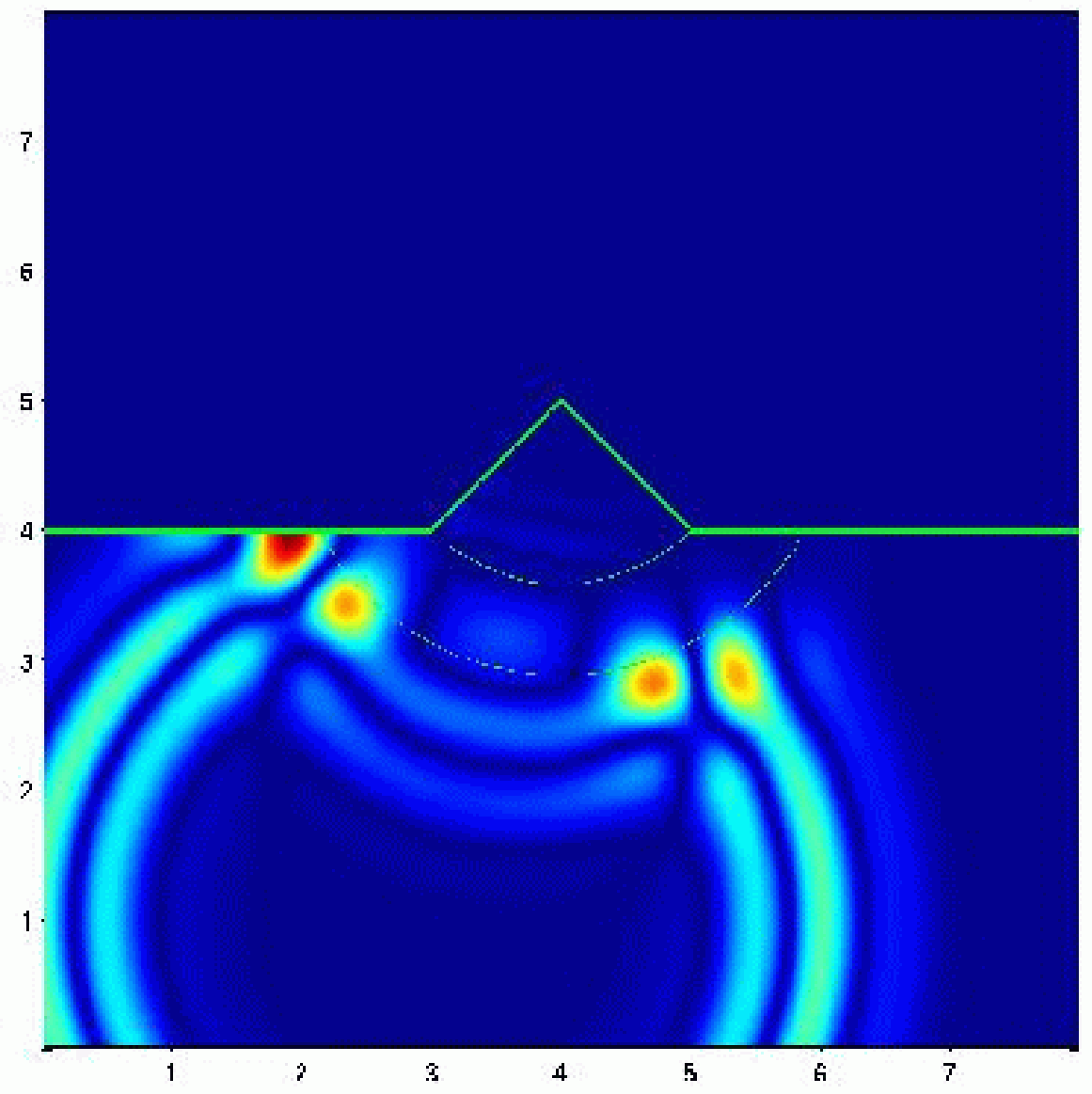,width=7cm}
\end{minipage}
\end{minipage}
~\\
\begin{minipage}{15cm}
\begin{minipage}{7cm}
\hspace*{3.2cm} $t=5$s
\end{minipage}
\begin{minipage}{7cm}
\hspace*{3.2cm} $t=8$s
\end{minipage}
\end{minipage}
~\\
\begin{minipage}{15cm}
\begin{minipage}{7cm}
\centering\psfig{figure=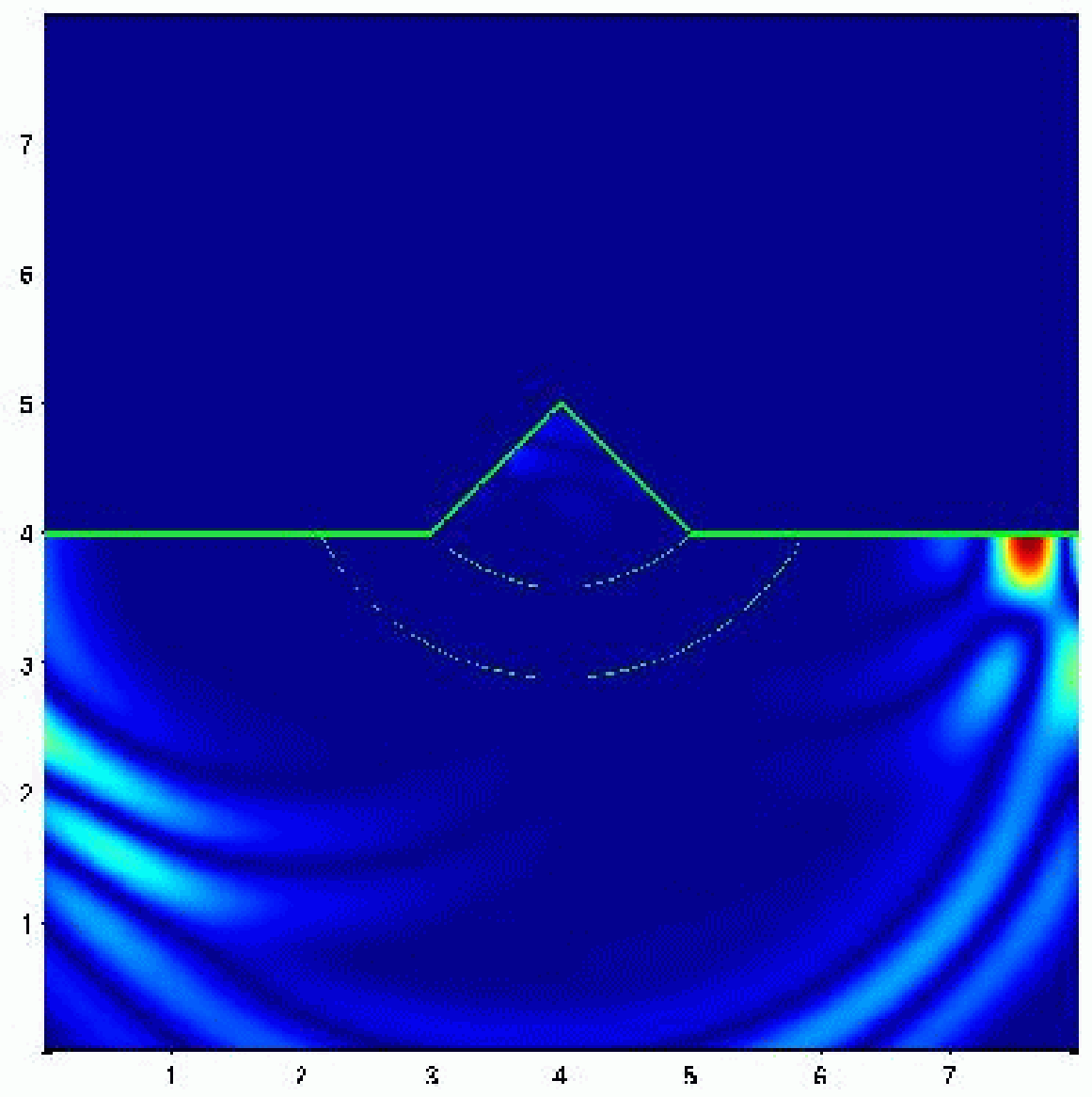,width=7cm}
\end{minipage}
\begin{minipage}{7cm}
\centering\psfig{figure=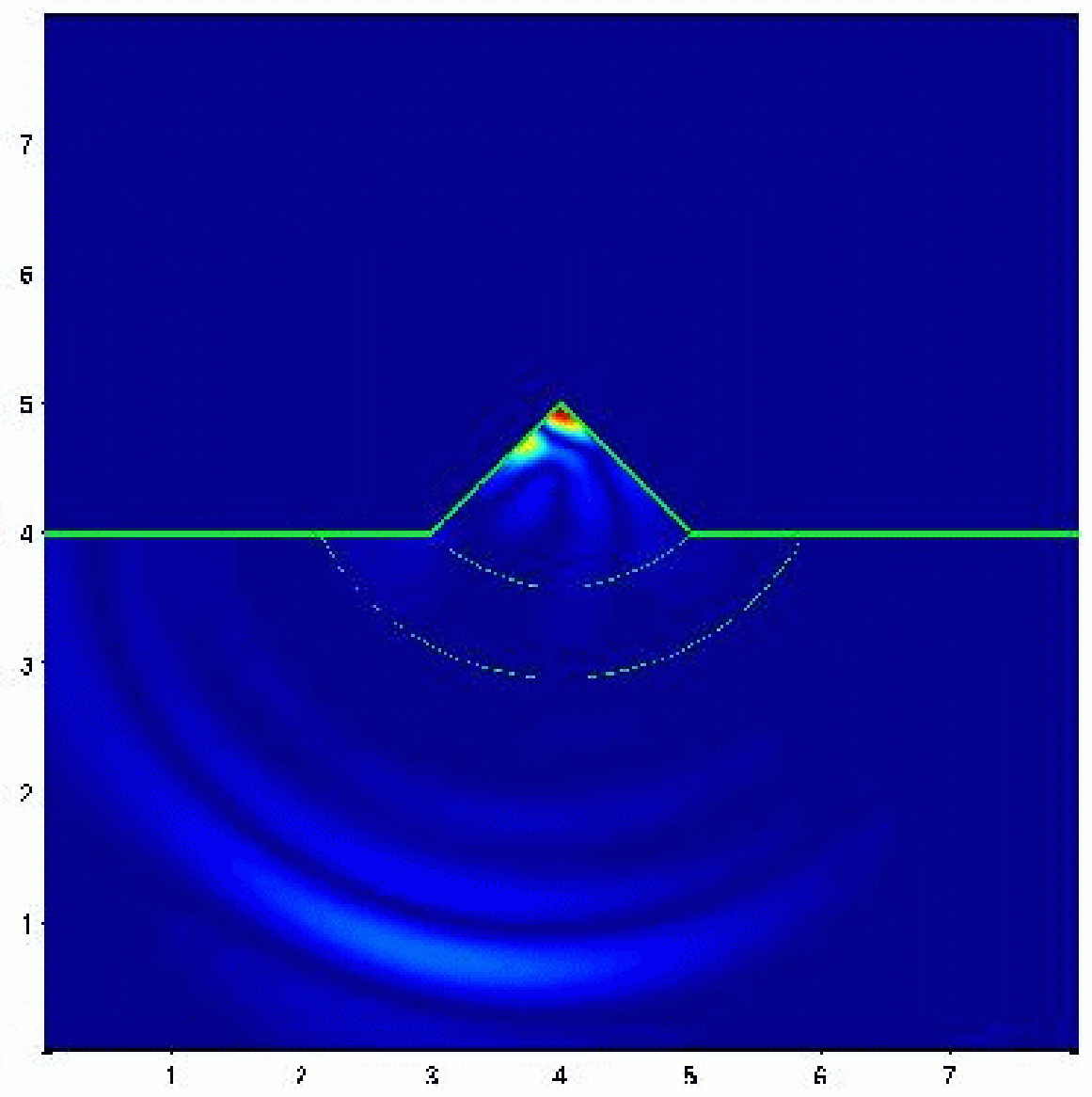,width=7cm}
\end{minipage}
\end{minipage}
~\\
\begin{minipage}{15cm}
\begin{minipage}{7cm}
\hspace*{3.2cm} $t=9$s
\end{minipage}
\begin{minipage}{7cm}
\hspace*{3.2cm} $t=9.5$s
\end{minipage}
\end{minipage}
~\\
\begin{minipage}{15cm}
\begin{minipage}{7cm}
\centering\psfig{figure=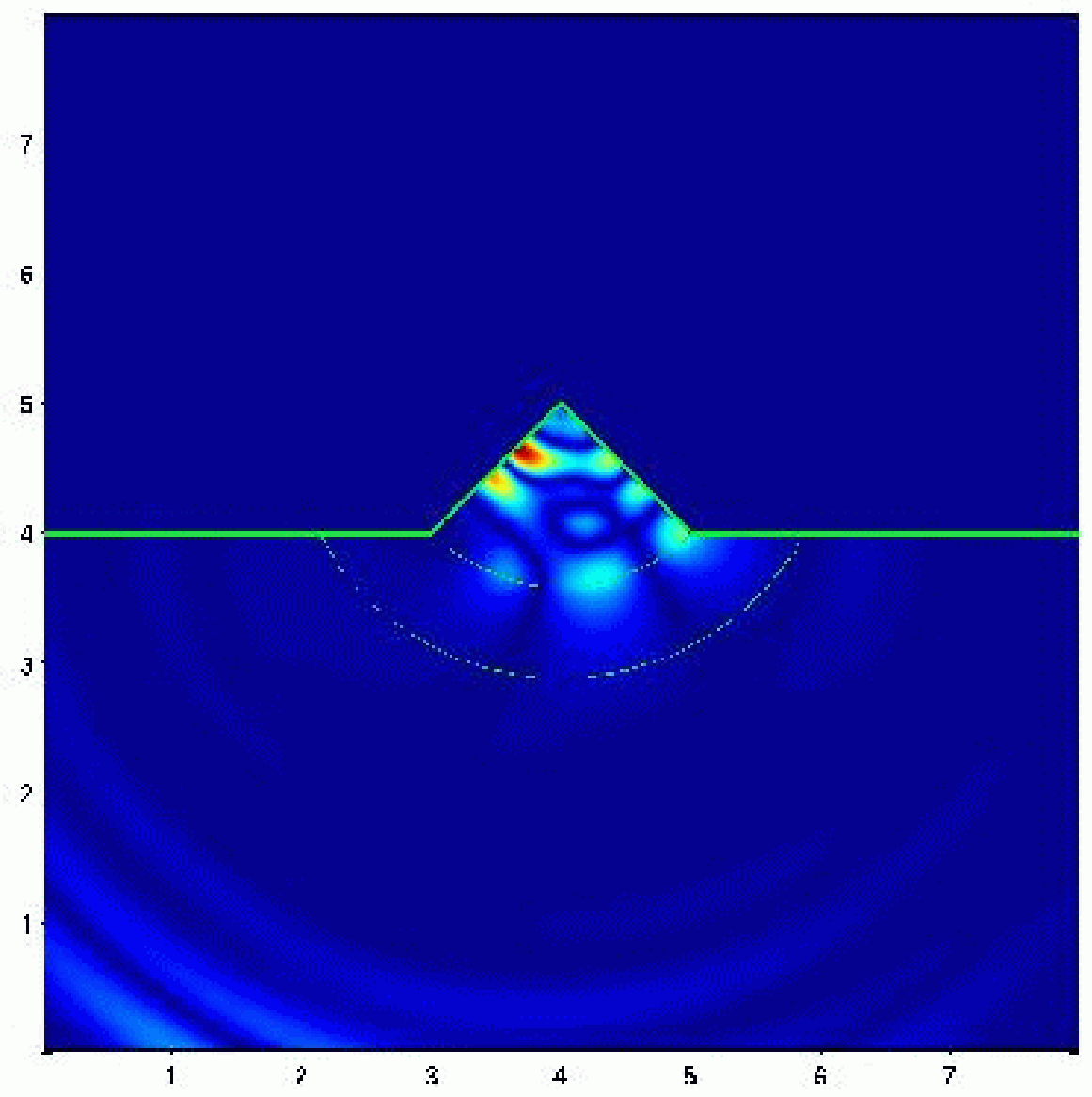,width=7cm}
\end{minipage}
\begin{minipage}{7cm}
\centering\psfig{figure=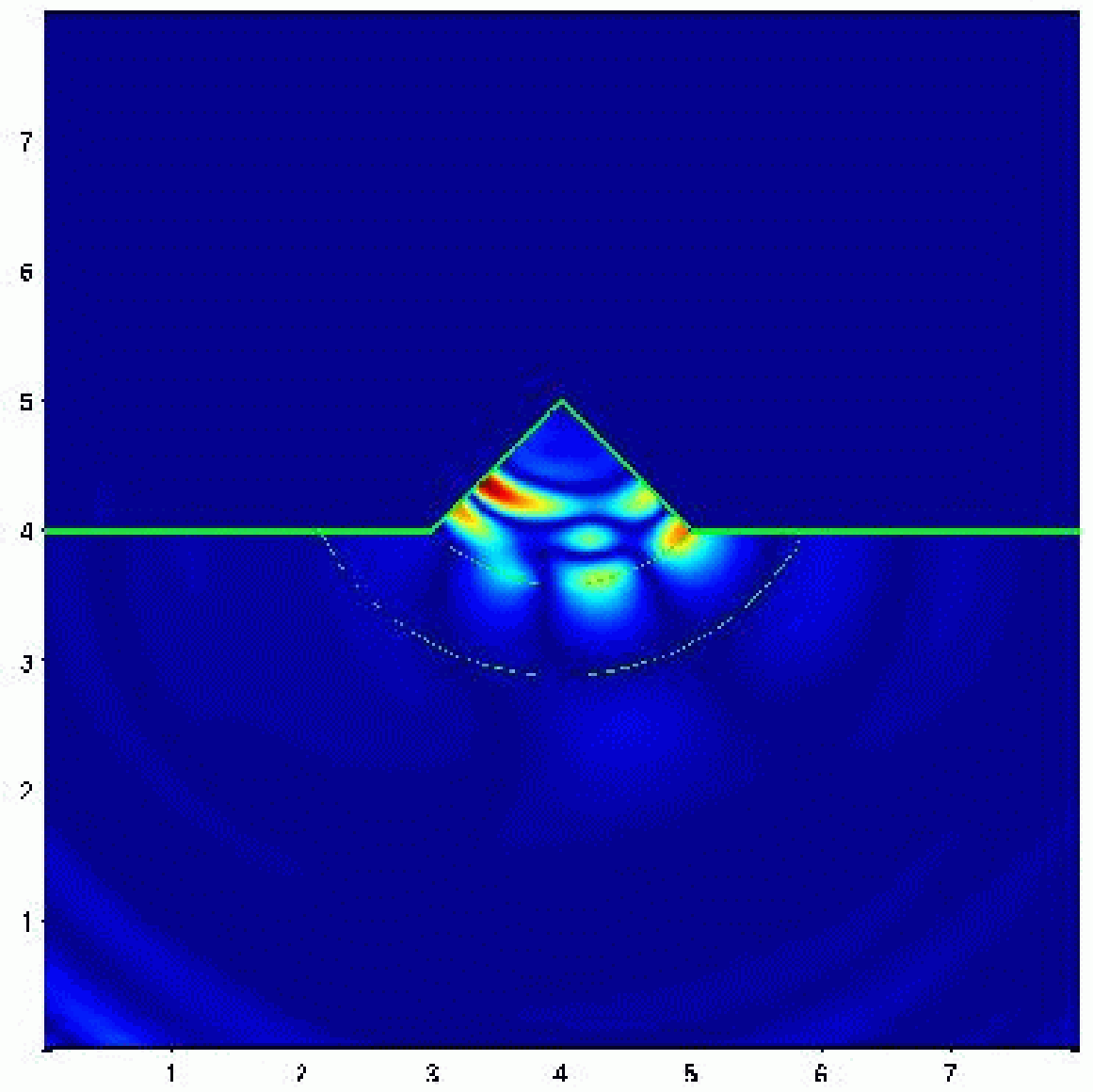,width=7cm}
\end{minipage}
\end{minipage}
\caption{Snapshots of the solution: the pressure field in the
  computational domain at different times}
\label{Montes}
\end{figure}

In Figure \ref{Montes} we display snapshots of the solution (the
pressure field) at different times. Diffraction from the free
surface is modeled by embedding the solution in a domain of a
simple shape using the fictitious domain method. To model wave
propagation in the infinite half-space the fictitious domain is
surrounded by an absorbing medium using the PML model. Although
for this problem a semi-analytical method similar to the one used
in \cite{Todo2001} can be employed to compute the solution, the
numerical method proposed in this paper is more general in that it
can be applied to any complicated geometry and/or propagation
media. Moreover, our numerical method can be of particular
interest in cases where the rigid foundation assumption is not
valid as it can provide realistic values for the stresses in the
structure.

\begin{center}
{\Large{Conclusion}}{\label{Conc}}
\end{center}
We employed a rational approximation of the frequency-dependent
viscoacoustic modulus in order to introduce dissipation into
time-domain computations. To do so, we followed the approach in
\cite{EK87} and chose relaxation frequencies $w_{l}(\mathbf{x})$
equidistant on a logarithmic scale in the frequency range
$[\frac{w_{max}}{100}; w_{max}]$, where $w_{max}$ is the maximal
frequency of the used source spectrum. This approach will be
accurate for propagation in media with a quality factor greater
than $10$. For media with high attenuation ($Q < 10$) it is
necessary in order to obtain accurate results to use a non-linear
minimization method such as the one proposed in \cite{Asv2003}.

By introducing this approximation of the viscoacoustic modulus
into the velocity-pressure formulation we obtained a first-order-
in-time linear system of equations. To discretize this system we
used a mixed finite-element method for the discretization in space
and a second-order finite difference scheme in time.

The velocity-pressure formulation was coupled with the fictitious
domain method in order to model the free surface boundary
condition on boundaries with complicated geometries, and with the
PML method to simulate wave propagation in unbounded domains. The
efficiency of the method was illustrated by numerical results.

\newpage \appendix
\section{Stability analysis}{\label{Ap1}
\subsection{The continuous problem}\label{Ap11}
We rewrite the continuous system in time with zero source term,
\begin{eqnarray}
\displaystyle \rho \frac{\partial \mathbf{v}}{\partial t}=\nabla
p, \label{Ap1Cont1}\\
\displaystyle \frac{\partial p}{\partial t} - \sum_{l=1}^{n}
\frac{\partial \eta_{l}}{\partial t} =
\mu_{R}\Div(\mathbf{v}), \label{Ap1Cont2} \\
\displaystyle \frac{\partial \eta_{l}}{\partial t} +
\omega_{l}\eta_{l} = \mu_{R}y_{l}\Div(\mathbf{v}), \forall l. \label{Ap1Cont3}
\end{eqnarray}
By taking the inner products (in $L^2$) of (\ref{Ap1Cont1}) with
$\mathbf{v}$,
(\ref{Ap1Cont2}) with $\dsp \left( p-\sum_{l=1}^{L}\eta_{l} \right)$,
and (\ref{Ap1Cont2}) with $\eta_{l}$ we get
\begin{eqnarray}
\displaystyle \left( \rho \frac{\partial \mathbf{v}}{\partial
  t},\mathbf{v}\right)=\left(\nabla (p),\mathbf{v}\right),
\label{Ap1Cont4}\\
\displaystyle \left(\frac{\partial}{\partial t}\left(p - \sum_{l=1}^{n}
 \eta_{l}\right),\left(p - \sum_{l=1}^{n}
 \eta_{l}\right)\right) = \mu_{R}\left( \Div(\mathbf{v}),\left(p -
 \sum_{l=1}^{n} \eta_{l}\right)\right), \label{Ap1Cont5} \\
\displaystyle \left(\frac{\partial \eta_{l}}{\partial t},\eta_{l}\right) +
\left(\omega_{l}\eta_{l},\eta_{l}\right) = \mu_{R}y_{l}\left(
 \Div(\mathbf{v}),\eta_{l}\right). \label{Ap1Cont6}
\end{eqnarray}
Then, summing $\dsp
(\ref{Ap1Cont4})+\frac{(\ref{Ap1Cont5})}{\mu_{R}}+\sum_{l}^{L}\frac{(\ref{Ap1Cont6})}{y_{l}\mu_{R}}$,
we obtain,
\begin{equation}
\left( \rho \frac{\partial v} {\partial t},\mathbf{v} \right) + \frac{1}{\mu_{R}} \left( \frac{\partial}{\partial t} \left( p - \sum_{l=1}^{n}
 \eta_{l}\right) ,\left( p - \sum_{l=1}^{n}
 \eta_{l} \right)  \right)+ \sum_{l=1}^{L} \frac{1}{y_{l}\mu_{R}} \left( \frac{\partial \eta_{l}}{\partial t},\eta_{l}\right)=-\sum_{l=1}^{L}\frac{\omega_{l}}{\mu_{R}y_{l}}\left(\eta_{l},\eta_{l}\right)
\label{Ap1Cont7}
\end{equation}
Keeping in mind that the energy of the system is,
\begin{equation}
\varepsilon = \frac{1}{2} \left( \rho \mathbf{v}, \mathbf{v} \right) +\frac{1}{2 \mu_{R}} \left( (p-\sum_{l=1}^{L}\eta_{l}) , (p-\sum_{l=1}^{L}\eta_{l}) \right) +\sum_{l=1}^{L} \frac{1}{2 y_{l} \mu_{R}} \left( \eta_{l},\eta_{l} \right)
\label{Ap1ECont}
\end{equation}
we finally get,
\begin{equation}
\frac{\partial \varepsilon}{\partial t} =
-\sum_{l=1}^{L}\frac{\omega_{l}}{\mu_{R}y_{l}}\|\eta_{l}\|^{2}\leq0 .
\end{equation}
Which implies that the energy of the system is decreasing with time,
when $\omega_{l}$, $\mu_{R}$ and $y_{l}$ are
positive quantities.  The relaxation frequencies $\omega_l$ are always
positive and the same holds for the relaxed modulus $\mu_R$. The
coefficients $y_l$ can in practice become negative if we do not solve
a constraint minimization problem. However,
we never encountered in practice a case for which
$$\dsp -\sum_{l=1}^{L} \frac{\omega_{l}}{\mu_{R}
  y_{l}}\|\eta_{l}\|^{2} \geq 0,$$
and thus the problem becomes unstable (in the sense that the
energy increases). To avoid this instability a constraint
minimization
  algorithm seeking for non-negative $y_l$ can be used.
\subsection{The discrete problem}{\label{Ap12}
We consider here the more general case where the pressure field is
discretized in $M_h^1$. Let us remark that the discretization
space $M_h^1$ admits the following orthogonal decomposition in
$L^2$, $$ M_h^1 = M_h \oplus (M_h)^{\perp},$$ where $M_h$ is the
space of piecewise constant functions, $$M_h = \left\{ {q}_h \in
L^2 \left/ \right. \forall  K \in  {\cal T}_{h},  {q}_{h\left|_{K}
    \right.} \in {P_0}(K) \right\},$$
and $(M_h)^{\perp}$ is its orthogonal complement in $M_h^1$ (with
respect to the inner product in $L^2$). To simplify the notation,
we denote by $P$ the discrete unknown associated with the pressure
field $P=p_h \in M_{h}^{1}$, so that we can write
$P=[P_{0},P_{1}]$ with $P_{0}$, the projection of $P$ on $M_h$ and
$P_{1}$ the projection of $P$ on $M_h^{\perp}$. The memory
variables are only discretized on $M_h$ .\\ In this case, we can
rewrite the discrete system as, (capital letters are used for the
discrete unknowns and the subscript $h$ is omitted)
\begin{eqnarray}
\displaystyle \rho \frac{V^{n+\frac{1}{2}}-V^{n-\frac{1}{2}}}{\Delta
  t}=-B^{0}_{h}P_{0}^{n}-B_{h}^{1}P_{1}^{n} \hspace{2.85cm}
\label{Ap1Dsi1}\\
\displaystyle \frac{P_{0}^{n+1}-P_{0}^{n}}{\Delta
  t}-\sum_{l=1}^{L}\frac{H_{l}^{n+1}-H_{l}^{n}}{\Delta
  t}=\mu_{R}B_{h}^{0,T}V^{n+\frac{1}{2}}\hspace{0.5cm}
\label{Ap1Dsi2} \\
\displaystyle \frac{P_{1}^{n+1}-P_{1}^{n}}{\Delta
  t}=\mu_{R}B_{h}^{1,T}V^{n+\frac{1}{2}}
\hspace{3.55cm}\label{Ap1Dsi3} \\
\displaystyle \frac{H_{l}^{n+1}-H_{l}^{n}}{\Delta t}+\omega_{l}
\frac{H_{l}^{n+1}+H_{l}^{n}}{2}= \mu_{R} y_{l}
B_{h}^{0,t}V^{n+\frac{1}{2}} \label{Ap1Dsi4}\hspace{0.7cm}
\end{eqnarray}
Then considering the inner products $( (\ref{Ap1Dsi1}) \text{ at time
  (n+1)} - ((\ref{Ap1Dsi1}) \text{ at time } n ) ) \times
  V^{n+\frac{1}{2}}$,\\ $\dsp (\ref{Ap1Dsi2}) \times
  \left(P_{0}^{n+1}-\sum_{l=1}^{L} H_{l}^{n+1}
  +P_{0}^{n}-\sum_{l=1}^{L} H_{l}^{n}  \right)$, $(\ref{Ap1Dsi3}) \times
  \left(P_{1}^{n+1}+P_{1}^{n} \right)$, and $(\ref{Ap1Dsi4}) \times
  \left(H_{l}^{n+1}+H_{l}^{n} \right)$, we get,
\begin{eqnarray}
\dsp \left(\rho V^{n+\frac{3}{2}},V^{n+\frac{1}{2}}
\right)=\hspace{13.3cm}\label{Ap1Dsi5}\\
\dsp \left(\rho V^{n+\frac{1}{2}},V^{n-\frac{1}{2}} \right)-\Delta t
\left(B_{h}^{0}(P_{0}^{n}+P_{0}^{n+1}),V^{n+\frac{1}{2}} \right)
-\Delta t \left(B_{h}^{1}(P_{1}^{n}+P_{1}^{n+1}),V^{n+\frac{1}{2}}
\right)\hspace{1.8cm} \nonumber\\
\dsp \|P_{0}^{n+1}-\sum_{l=1}^{L}H_{l}^{n+1}
\|^{2}=\hspace{12.5cm}\label{Ap1Dsi6}\\
\dsp \|P_{0}^{n}-\sum_{l=1}^{L}H_{l}^{n} \|^{2} + \Delta t \ \mu_{R}
\left(B_{h}^{0,T}V^{n+\frac{1}{2}},P_{0}^{n+1}+P_{0}^{n}
\right)-\Delta t \mu_{R}
\left(B_{h}^{0,T}V^{n+\frac{1}{2}},\sum_{l=1}^{L}H_{l}^{n+1}+H_{l}^{n}
\right)\nonumber\\
\dsp \|P_{1}^{n+1} \|^{2}=\|P_{1}^{n} \|^{2}+\Delta t \mu_{R} \left(
B_{h}^{1,T} V^{n+\frac{1}{2}},P_{1}^{n+1} +P_{1}^{n}
\right)\hspace{7.5cm}\label{Ap1Dsi7}\\
\dsp \|H_{l}^{n+1} \|^{2}=\|H_{l}^{n} \|^{2}-\omega_{l} \Delta t
\frac{ \|H_{l}^{n+1}-H_{l}^{n}\|^{2}}{2}+ \mu_{R} \Delta t y_{l}\left(
B_{h}^{0,T}V^{n+\frac{1}{2}},H_{l}^{n+1}+H_{l}^{n}
\right)\hspace{3.35cm}\label{Ap1Dsi8}
\end{eqnarray}
Finally summing $\dsp (\ref{Ap1Dsi5})+\frac{(\ref{Ap1Dsi6})}{\mu_{R}}+
\frac{(\ref{Ap1Dsi7})}{\mu_{R}}
+\sum_{l=1}^{L}\frac{(\ref{Ap1Dsi8})}{y_{l}\mu_{R}}$, we get,
\begin{equation}
\frac{\varepsilon_{h}^{n+1}-\varepsilon_{h}^{n}}{\Delta t }= -
\sum_{l}^{L}\frac{\omega_{l}}{\mu_{R}
  y_{l}}\frac{\|H_{l}^{n+1}+H_{l}^{n}\|^{2}}{4},
\label{en.d}
\end{equation}
with the discrete energy being defined by,
\begin{equation}
2\varepsilon_{h}^{n}=\left( \rho
V^{n+\frac{1}{2}},V^{n-\frac{1}{2}}\right)+
\frac{1}{\mu_{R}}\|P_{1}^{n}\|^{2}+
\frac{1}{\mu_{R}}\|P_{0}^{n}-
\sum_{l=1}^{L}H_{l}^{n}\|^{2}+\sum_{l=1}^{L}
\frac{1}{\mu_{R}y_{l}}\|H_{l}^{n}\|^{2}.
\label{Endis}
\end{equation}
Equation (\ref{en.d}) shows that the discrete energy is also
decreasing, under the same assumptions on $y_{l}$ as in \ref{Ap11}.

To show under which condition the quantity defined by (\ref{Endis}) is
positive and thus an energy, we use the orthogonality relation between
$P_{0}$ and $P_{1}$ (note that $P_1$ is also orthogonal to $H_l$), to
get,
$$\begin{array}{l}
\dsp 2\varepsilon_{h}^{n}=\left( \rho
(V^{n+\frac{1}{2}}+V^{n-\frac{1}{2}}),(V^{n+\frac{1}{2}}
+V^{n-\frac{1}{2}})\right)+\frac{1}{\mu_{R}}\|P^{n}\|^{2}
+\sum_{l=1}^{L}\frac{1}{\mu_{R}y_{l}}\|H_{l}^{n}\|^{2}+\\
\dsp \frac{1}{\mu_{R}}\left( \sum_{l=1}^{L}
H_{l}^{n},\sum_{l=1}^{L}H_{l}^{n}
\right)-\frac{2}{\mu_{R}}\left(P^{n},\sum_{l=1}^{L}H_{l}^{n}
\right)-\frac{\Delta t^{2}}{4\rho}\left(B_{h}P^{n}, B_{h}P^{n}\right)
\end{array}
$$
or
$$
\dsp 2\varepsilon_{h}^{n} \geq \frac{1}{\mu_{R}} \left[
  \left(1-\frac{\Delta t^{2}\mu_{R} \| B_h \|^2 }{4\rho}
  \right)\|P^{n}\|^{2}+ \left\| \sum_{l=1}^{L}
  H_{l}^{n} \right\| - 2 \left(P^{n},\sum_{l=1}^{L}H_{l}^{n}
\right)
+
  \sum_{l=1}^{L}\frac{1}{y_{l}}\|H_{l}^{n}\|^{2}\right]
$$ where $P^{n}=P_{0}^{n}+P_{1}^{n}$ and
$B_{h}P^{n}=B_{h}^{0}P_{0}^{n}+B_{h}^{1}P_{1}^{n}$. We rewrite
this equation as a matrix associated with the quadratic
formulation and we prove that the eigenvalues of this matrix are
positive under the CFL condition,
\begin{equation}
\frac{\Delta t^{2}}{4} \frac{\mu_{R}}{\rho} |\!| B_{h} |\!| ^{2}
\left( 1+ \sum_{l=1}^{L} y_{l} \right) \leqq 1
\end{equation}
with $ |\!| B_{h}^{T} B_{h} |\!| \geqq \frac{4}{h^{2}}$ in 1D and $
|\!| B_{h}^{T} B_{h} |\!| \geqq \frac{8}{h^{2}}$ in 2D.
\section{Dispersion analysis}{\label{Ap2}}
\subsection{The continuous problem}{\label{Ap21}}
Suppose that $\mathbf{v}(\mathbf{x},t)$, $p(\mathbf{x},t)$, and
$\eta_{l}(\mathbf{x},t)\forall l$, are plane waves,
$$\left\{
\begin{array}{ll}
\dsp \mathbf{v}(\mathbf{x},t)=\mathbf{v_{0}}{\exp{\left( \mathbf{i}
    \left(\omega t - \mathbf{K}\mathbf{x} \right) \right)}}, \\
\dsp p(\mathbf{x},t)=p_{0}{\exp{\left( \mathbf{i} \left(\omega t -
    \mathbf{K}\mathbf{x} \right) \right)}}, \\
\dsp \eta_{l}(\mathbf{x},t)=\eta_{l}^{0}{\exp{\left( \mathbf{i}
    \left(\omega t - \mathbf{K}\mathbf{x} \right) \right)}},
\end{array}
\right.
$$
where $\mathbf{K}\mathbf{x}=kx$ in 1D and $\mathbf{K}
\mathbf{x}=k_{x}x+k_{y}y=k{\cos(\Phi)}x+k{\sin(\Phi)}y$,  $\Phi$ being
the incident angle of the plane wave in 2D. Introducing this
expression into the time domain system (\ref{systime}), we get the
dispersion relation,
\begin{equation}
\omega^{2}=\mathbf{K}^{2}c_{R}^{2}\left(1+\sum_{l=1}^{L}\frac{\mathbf{i}\omega
  y_{l}}{\mathbf{i}\omega+\omega_{l}} \right)
\label{DispContAp}
\end{equation}
with $c_{R}=\sqrt{\frac{\mu_{R}}{\rho}}$ the relaxed velocity. If
the medium is non-dissipative (i.e., $y_{l}=0 \forall l$),
(\ref{DispContAp}) becomes the well-known relation
$\omega^{2}=\mathbf{K}^{2}c^{2}$. Note that the dispersion
relation (\ref{DispContAp}) is no longer explicit in $\omega$.
\subsection{The discrete problem}{\label{Ap22}}
We are interested in the general formulation for which the
pressure field is discretized in $M_{h}^{1}$ and $\eta_{l}$ in
$M_{h}$. Considering that $V$, $P$, and $H_{l}$ are plane waves,
and employing the same notation as in \ref{Ap12}, we get,
\begin{equation}
\dsp {\sin^{2}\left(\chi_{t} \right)} =\frac{\Delta t^{2}
  c_{R}^{2}}{4}\left(B_{h}B_{h}^{T}+\sum_{l=1}^{L} B_{h}^{0}
B_{h}^{0,T} \frac{2\mathbf{i} y_{l} {\tan\left(
    \chi_{t}\right)}}{\Delta t
  \omega_{l}+2\mathbf{i}{\tan\left(\chi_{t } \right)}}\right)
\end{equation}
wherein $\chi_{t}=\frac{\omega \Delta t}{2}$, $\Delta t$ being the
discretization step in time. After some calculations we obtain,
\begin{equation}
\begin{array}{ll}
\dsp \sin^{2} \left( \frac{\omega \Delta_{t}}{2}  \right)=\frac{\Delta_{t}^2 c^2}{4} \left({\sin^{2} \left(\frac{k_{x}\Delta_{x}}{2} \right)}+{\sin^{2}\left(\frac{k_{y}\Delta_{y}}{2} \right)} \right) \left(1+\sum_{l=1}^{L} \frac{2\mathbf{i} y_{l} \tan \left(\frac{\omega \Delta_{t}}{2}  \right)}{\Delta_{t}\omega_{l}+2 \mathbf{i} \tan \left(\frac{\omega \Delta_{t}}{2}  \right) } \right) & \dsp \text{ in 2D}\\
\dsp \sin^{2} \left( \frac{\omega \Delta_{t}}{2}  \right)=\frac{\Delta_{t}^2 c^2}{4}\left({\sin^{2} \left(\frac{k\Delta_{x}}{2} \right)}\right) \left(1+\sum_{l=1}^{L} \frac{2\mathbf{i} y_{l} \tan \left(\frac{\omega \Delta_{t}}{2}  \right)}{\Delta_{t}\omega_{l}+2 \mathbf{i} \tan \left(\frac{\omega \Delta_{t}}{2}  \right) } \right) & \dsp \text{ in 1D}
\end{array}
\end{equation}
wherein $\Delta_{x}$ and $\Delta_{y}$ are the discretization step
in space. In our case $\Delta_{x}=\Delta_{y}=h$.

\bibliography{biblio}
\bibliographystyle{plain}

\end{document}